\documentclass{article}

\usepackage[T1]{fontenc}

\usepackage[a4paper,left=2cm,right=2cm,top=3cm,bottom=3cm]{geometry}
\usepackage{csquotes}
\usepackage{orcidlink}

\usepackage{amsmath,amsthm,amssymb}

\usepackage{nicefrac}

\usepackage{hyperref}

\usepackage{tikz}
\usepackage{pgfplots}
\usetikzlibrary{pgfplots.groupplots}
\pgfplotsset{compat=1.17}
\usepackage{pgfplotstable}

\usepackage{placeins}

\usepackage
[
style=iso-alphabetic,
maxnames=5, 
giveninits=true, 
isbn=false,url=false,eprint=false, doi=true,
backref=true
]
{biblatex}

\DeclareSourcemap{
  \maps[datatype=bibtex]{
    \map[overwrite]{
      \step[fieldsource=doi, final]
      \step[fieldset=url, null]
      \step[fieldset=eprint, null]
    }  
  }
}

\NewBibliographyString{toappearin}
\NewBibliographyString{submittedto}
\DefineBibliographyStrings{english}{%
  toappearin  = {to appear in},
  submittedto = {submitted to},
}

\renewbibmacro*{in:}{%
  \ifboolexpr{not test {\iffieldundef{pubstate}}
              and (test {\iffieldequalstr{pubstate}{toappearin}}
                   or test{\iffieldequalstr{pubstate}{submittedto}})}
    {\printtext{%
       \printfield{pubstate}\intitlepunct}%
     \clearfield{pubstate}}
    {\printtext{%
       \bibstring{in}\intitlepunct}}}

\usepackage[final]{showkeys} 
\usepackage{todonotes}

\theoremstyle{plain}
\newtheorem{theorem}{Theorem}[section]
\newtheorem{lemma}[theorem]{Lemma}
\newtheorem{corollary}[theorem]{Corollary}
\newtheorem{assumption}[theorem]{Assumption}

\theoremstyle{definition}
\newtheorem{definition}[theorem]{Definition}

\theoremstyle{remark}
\newtheorem{remark}[theorem]{Remark}

\DeclareMathOperator{\Span}{span}
\DeclareMathOperator{\sign}{sign}

\newcommand{\RR}{\mathbb{R}}

\newcommand{\dx}{\,\mathrm{dx} }

\allowdisplaybreaks

\begin{document}

\title{A least squares space-time approach for parabolic equations}

\author{
Michael Hinze%
\footnote{Mathematical Institute, University of Koblenz, Germany}
\orcidlink{0000-0001-9688-0150},\,%
Christian Kahle%
\footnotemark[1]
\orcidlink{0000-0002-3514-5512},\,%
Michael Stahl%
\footnotemark[1]
\orcidlink{0009-0000-4692-7328}
}

\date{\today}

\maketitle

\begin{abstract}
We propose a least squares formulation for abstract parabolic equations in the natural $L^2(0,T;V^\star)\times H$ norm which only relies on natural regularity assumptions on the data of the problem. The resulting bilinear form then is symmetric, coercive and continuous. We provide two space-time Galerkin frameworks for the numerical approximation. The first one uses a conformal discretization of the underlying bilinear system and relies on the fact that the $V^*-$norm of basis functions can be evaluated exactly. The second approach is nonconforming an replaces the evaluation of the $V^*-$norm by a discrete pendant. We prove convergence for both approaches and illustrate our analytical findings by selected numerical experiments.
\end{abstract}

\section{Introduction.}

Least squares formulations of partial differential equations (PDEs) have been receiving more attention for some time now. In the context of the numerical treatment of the heat equation within a holistic space-time framework, the least squares approach for data in $L^2(0,T;L^2(\Omega))\times L^2(\Omega)$ is introduced in \cite{FuehrerKarkulik-SpaceTimeparabolic} by formulating the elliptic part in a mixed form. This enables work with a variational formulation whose associated bilinear form is coercive in space-time, thus allowing the standard Galerkin framework to be applied.

\paragraph{Contribution:} In this study, we examine the least squares formulation of parabolic problems in the natural $L^2(0,T;H^{-1}(\Omega))\times L^2(\Omega)$ norm, which only relies on natural regularity assumptions on the data. The resulting bilinear form then is symmetric, continuous and coercive. We propose a space-time Galerkin approach for this formulation, which yields a conforming discretization with quasi-optimal approximation properties. For its practical numerical implementation, we restrict to one spatial dimension of the domain. In this case, the negative Sobolev norms that appear in the bilinear form are evaluated analytically. We present numerical results that confirm the quasi-optimal approximation properties of this approach. In addition, we propose a nonconforming discretization where dual norms are approximately evaluated using an approximation of the Riesz map. For this approach we prove convergence and develop tailored preconditioners for the resulting discrete space-time systems.

\medskip
\paragraph{Novelty statement:}
\begin{itemize}
    \item[$\blacktriangleright$] Least squares space-time formulation of parabolic problems requiring only minimal regularity (Section \ref{sec:setting});
    \item[$\blacktriangleright$] Quasi-optimal approximation property of the conformal Galerkin approximation with numerical implementation in one spatial dimension (Section~\ref{GA});
    \item[$\blacktriangleright$] Convergence of a nonconforming Galerkin approximation relying on a saddle point formulation (Section \ref{sec:FESP});
    \item[$\blacktriangleright$] Development of quasi-optimal preconditioners for the saddle point formulation (Section \ref{sec:FE:num}).
\end{itemize}

\paragraph{Literature:} Many contributions have been made to a space-time approach for parabolic PDEs.
We focus on Petrov--Galerkin and least squares schemes, and refer to \cite{andreev2013stability,Bochev-2014-minres-ls-FEM,gander2016analysis,langer2019space,langer2021efficient,perugia2023exponential,schwab2009space,stevenson2021stability,steinbach2015space}.
Least squares formulations for elliptic PDEs are often based on the FOSLS approach.
This assumes a higher regularity of the solution, which allows a mixed form to be used for the elliptic operator. \cite{BochevGunzburger-LeastSquares}. 
This approach is typically not used in combination with natural regularity of data. 
In \cite{Bramble1997} the problem of assuming higher regularity is addressed by considering a first-order system again, but this time using the natural $H^{-1}$-norm of the objects under consideration in the least squares formulation. The authors also provide an efficient preconditioning technique for evaluating the norm. They also discuss error estimates.
Another example, where natural regularity is employed, is provided by \cite{Fuehrer2022,Fuehrer2024-MixedFEMFOSLS-Hm1-loads}.
In this case, a regularization technique is employed to address volume forces with natural $H^{-1}$-regularity, resulting in optimal convergence rates for both a Petrov–Galerkin scheme and a FOSLS scheme.
In \cite{FuehrerKarkulik-SpaceTimeparabolic} a full space-time least squares approach for parabolic equations is proposed. To achieve this, the FOSLS approach is extended directly to time-dependent problems by incorporating time as an additional dimension in the least squares problem. In this work, the authors demonstrate quasi-optimality of the numerical approximation using space-time finite elements, and provide approximation estimates for an approximation using piecewise globally continuous linear elements. They present numerical results demonstrating that this formulation provides full space-time adaptivity.
The relevant assumption for \cite{FuehrerKarkulik-SpaceTimeparabolic} is that the elliptic operator has higher regularity in order to obtain the mixed formulation. 
We also refer to \cite{Steinbach2023-adaptiveleastsquaresspacetimefinite} for a recent alternative approach to consider general PDEs in the least squares sense and in natural norms. In an abstract setting, the natural norm is realised by a suitable operator, and minimisation leads to a saddle point problem that includes the PDE operator and its adjoint. This approach naturally results in an indicator for the space-time adaptivity of the discretization mesh.

Of course are there also contributions to space-time formulations of other than parabolic equations, see e.g.~\cite{dorfler2020space} for hyperbolic systems, \cite{hain2024ultra} for a ultra-weak formulation of the Schrödinger equation, and \cite{henning2022ultraweak,loscher2023numerical,steinbach2022generalized} for the wave equation.

\medskip

Our approach is useful in the context of reduced basis methods for parabolic PDEs, since it enables the direct  transfer of numerical techniques from the well-established and effective treatment of e.g.~parametrized elliptic PDEs to the parabolic setting. However, a posteriori error bounds for parabolic equations are often sub-optimal due to the equations being treated in a separate way in terms of space and time during the discretization process. Our approach enables classical residual error estimators to be used analogously to the elliptic setting, allowing certification in the natural $W(0, T)$-norm.
For similar reasons, our approach is also useful for the numerical treatment of optimal control problems involving parabolic PDEs, since the entire error analysis can be adapted from the well-developed elliptic case. Furthermore, it is particularly well suited to shape optimisation problems involving parabolic equations. In this context, the formulation allows parabolic PDE-constrained shape optimisation problems involving time-dependent functionals and shapes to be solved using the same techniques that have been developed for elliptic PDE-constrained shape optimisation. Related results on these topics will be published elsewhere.

\paragraph{Outline.} In Section \ref{sec:setting} we introduce the problem setting and derive a symmetric, coercive and continuous bilinear form based on the least squares problem.
We also introduce an equivalent saddle point formulation.
In Section~\ref{sec:1d} we investigate the one-dimensional case based on a conforming finite element discretization introduced in Section~\ref{sec:GA-conforming}. We then in Section \ref{sec:FESP} introduce a discretization of the  equivalent saddle point problem for which we propose  a convergent finite element approximation and also develop  preconditioners for the resulting discrete space-time systems. The performance of the method is investigated in Section \ref{sec:FE:num}.


\section{Problem setting.}
\label{sec:setting}

Let $(V,(\cdot,\cdot)_{V})$ and $(H,(\cdot,\cdot)_H)$ 
denote separable Hilbert spaces with the properties 
$V \hookrightarrow H\equiv H^\star \hookrightarrow V^\star$, so that $(V,H,V^\star)$ forms a Gelfand triple, and $a : V \times V \rightarrow \RR$ be a symmetric, continuous and coercive bilinear form, defining the inner product in $V$,
\begin{align}\label{eq:RBbilinA}
    (u, v)_V := a(u, v) \qquad \forall u, v \in V .
\end{align}
Then  $a(u,\cdot) \in V^\star$ denotes a linear form on $V$ for any $u \in V$.
With $A: V \rightarrow V^\star$ we denote the Riesz isomorphism, defined by 
\begin{align}
    \langle A u, v \rangle_{V^\star, V} := a(u, v) \qquad \forall u, v \in V.
\end{align}
From here onwards we call $R := A^{-1} : V^\star \rightarrow V$ the Riesz lift, satisfying for any $\phi \in V^\star$
\begin{align}
    (R \phi,v)_{V} =  \langle \phi,v\rangle_{V^\star,V}  \qquad \forall  v \in V
\end{align}
and set for $\phi,\psi \in V^\star$
\begin{align}
    (\phi,\psi)_{V^\star}:= (R \phi,R \psi)_{V}. 
\end{align}
By definition of $(\cdot,\cdot)_V$ it follows that $A R \phi = \phi$ for all $\phi\in V^\star$ and 
$R A v = v$ for all $v \in V $ and
since $(V,H,V^\star)$ forms a Gelfand triple the following equality holds
for any $ \phi\in H, \psi \in V^\star$
\begin{align}
\label{eq:RieszProp}
    (\phi,\psi)_{V^\star} = (R\phi,R\psi)_V = \langle \phi,R\psi\rangle_{V^\star,V} = (\phi,R\psi)_H,
\end{align}
see \cite{Wloka1987,Troeltzsch2010}.

\bigskip

For given $T>0$ we define the space $W(0,T):=\{v\in L^2(0,T;V)\mid v_t\in L^2(0,T;V^\star)\}$. 
The canonical inner product in $W(0,T)$ is $(u,v)_{W(0,T)}:=\int_0^T (u_t,v_t)_{V^\star} + \int_0^T(u,v)_{V}$ and the norm induced by this inner product is denoted by $\|v\|^2_{W(0,T)} = (v,v)_{W(0,T)}$. 

For given $f\in L^2(0,T;V^\star)$ and $y_0\in H$ we consider the parabolic equation
\begin{equation}\label{eq:parab}
    y_t + Ay = f \text{ in } L^2(0,T;V^\star), \quad y(0)=y_0 \text{ in } H.
\end{equation}
It is well known that \eqref{eq:parab} admits a unique variational solution $y\in W(0,T)$, see e.g.~\cite{Wloka1987,Thomee2006}.
In this context we naturally extent $R$ from $V^\star$ to $L^2(0,T;V^\star)$ and understand $R$ as operator $R:L^2(0,T;V^\star) \to L^2(0,T;V)$ and define $R\phi$ for any $\phi \in L^2(0,T;V^\star)$  as $R\phi(t)$ for almost every $t$. With this $R\phi$ fulfills
\begin{equation}
    \int_0^T (R\phi, v)_V = \int_0^T \langle \phi, v \rangle_{V^\star, V} \qquad \forall v \in V .
\end{equation}

\begin{remark}
The classical application for our setting is the heat equation with homogeneous Dirichlet boundary data. In this application for an open and bounded domain $\Omega \subset \mathbb{R}^n$ we set
$A := -\Delta$, $V:= H^1_0(\Omega)$, $H=L^2(\Omega)$ and $V^\star=H^{-1}(\Omega)$.
Then $\langle A y,v\rangle_{V^\star,V} = \int_\Omega \nabla y \nabla v $,
and the inner product on $V$ is given by  $(y,v)_V = \int_\Omega \nabla y\nabla v$.
\end{remark}

In the following we introduce a space-time least squares formulation of \eqref{eq:parab}. 

\begin{lemma}
Consider the minimization problem 
\begin{equation}
\tag{\ensuremath{\mathcal M_r}}
\label{eq:minresid}
    \min_{v\in W(0,T)} \|r(v)\|^2_{{ L^2(0,T;V^\star)}\times H} 
    := \|v_t+Av-f\|^2_{ L^2(0,T;V^\star)}+\|v(0)-y_0\|^2_H.
\end{equation}
Let $y \in W(0,T)$    denote the unique solution to \eqref{eq:parab}.
Then $y$ also is the unique solution to \eqref{eq:minresid} and satisfies
the necessary and sufficient first order optimality condition
\begin{align}
    \int_0^T ( y_t+Ay,w_t+Aw)_{V^\star} + (y(0),w(0))_H =
\int_0^T ( f,w_t+Aw)_{V^\star} + (y_0,w(0))_H\quad \forall w \in W(0,T).
\label{eq:beql-V1}
\end{align}
\end{lemma}
\begin{proof}
Let $y\in W(0,T)$ denote the unique solution to \eqref{eq:parab}.
Because $\|r(y)\|_{L^2(0,T;V^{\star})\times H} = 0$, problem \eqref{eq:minresid} admits at least one solution.
Since $y$ is unique and by construction, the solution to \eqref{eq:minresid} is unique.

The necessary and sufficient first order optimality condition \eqref{eq:beql-V1} follows from direct calculation. 
\end{proof}

Next, we reformulate \eqref{eq:beql-V1}.
To begin with, let $v,w \in W(0,T)$ be arbitrary.
Using the identities  
\begin{align}
   & \int_0^T \langle v_t,w\rangle_{V^\star,V} +  \int_0^T\langle w_t,v\rangle_{V^\star,V}
= (v(T),w(T))_H - (v(0),w(0))_H, \quad \quad
(v_t,Aw)_{V^\star} = \langle v_t, w\rangle_{V^\star,V},
\end{align}
a short calculation shows
\begin{align}
     &\int_0^T ( v_t+Av,w_t+Aw)_{V^\star} + (v(0),w(0))_H\\
    &= \int_0^T (v_t,w_t)_{V^\star} + \int_0^T(Av,Aw)_{V^\star}
     + \int_0^T(v_t,Aw)_{V^\star} 
     + \int_0^T(Av,w_t)_{V^\star}   + (v(0),w(0))_H\\
     &=\int_0^T (v_t,w_t)_{V^\star}+ \int_0^T(v,w)_{V}  + (v(T),w(T))_H.
\end{align}
Moreover it holds
\begin{align}
   & (y_0,w(0))_H + \int_0^T(f,w_t+Aw)_{V^\star}
    = (y_0,w(0))_H + \int_0^T ( f,w_t)_{V^\star}
    + \int_0^T \langle f,w\rangle_{V^\star,V}. 
\end{align}

\begin{definition}
We consider the following formulation to solve \eqref{eq:parab}:\\
Find $y \in W(0,T)$ such that
\begin{align} 
    \tag{\ensuremath{P}}
    b(y,w) = l(w) \quad \forall w \in W(0,T),
    \label{prob:P}
\end{align}
where for arbitrary $v,w\in W(0,T)$ the forms $b:W(0,T) \times W(0,T) \to \RR$ and $l:W(0,T)\to \RR$ are defined as
\begin{align}
  b(v,w) &:= \int_0^T (v_t,w_t)_{V^\star}+ \int_0^T(v,w)_{V}  + (v(T),w(T))_H, \label{eq:b_def}\\
    l(w)&:=(y_0,w(0))_H 
    + \int_0^T ( f,w_t)_{V^\star}
    + \int_0^T \langle f,w\rangle_{V^\star,V}. 
    \label{eq:l_def}
\end{align}
\end{definition}

\begin{theorem}
\label{thm:cont:bl-LaxMilgram}
The bilinear form $b:W(0,T) \times W(0,T) \to \RR$ is symmetric, bounded and coercive. Moreover the linear form $l:W(0,T)\to \RR$ is continuous.
\end{theorem}
\begin{proof}
For any $v,w \in W(0,T)$ it holds $b(v,w) = (v,w)_{W(0,T)} + (v(T),w(T))_H$. 
Now symmetry is obviously given.
To show coercivity, let $y\in W(0,T)$ be given. 
\begin{align}
    b(y,y) =(y,y)_{W(0,T)} + (y(T),y(T))_H \geq \|y\|^2_{W(0,T)}.
\end{align}

To show boundedness, let $y,w \in W(0,T)$ be arbitrary. Then with a constant $C>0$ that does not depend on $y,w$ it holds
\begin{align}
   | b(y,w)| &= (y,w)_{W(0,T)} + (y(T),w(T))_H \\
    &\leq \|y\|_{W(0,T)} \|w\|_{W(0,T)} + \|y(T)\|_H\|w(T)\|_H \\
    &\leq \|y\|_{W(0,T)} \|w\|_{W(0,T)} + \|y\|_{C([0,T];H)}\|w\|_{C([0,T];H)} \\
    &\leq \|y\|_{W(0,T)} \|w\|_{W(0,T)} + C\|y\|_{W(0,T)}\|w\|_{W(0,T)} \\
    &= (1+C)\|y\|_{W(0,T)} \|w\|_{W(0,T)},
\end{align}
where we used the continuous embedding $C([0,T],H) \hookrightarrow W(0,T)$.

The continuity of $l$ follows by similar arguments.
\end{proof}

\begin{remark}
With Theorem~\ref{thm:cont:bl-LaxMilgram} the existence of a unique solution $y\in W(0,T)$ to \eqref{prob:P} also directly follows from Lax--Milgram's theorem. 
\end{remark}

\subsection{Reformulation as saddle point problem.}

For the fully practical numerical analysis of problem \eqref{prob:P} presented in Section~\ref{sec:FESP} it is convenient to introduce an equivalent saddle point formulation. To begin with we define 
$\hat{a}:W(0,T)\times W(0,T) \to \RR$, 
$\hat{b}:W(0,T)\times L^2(0,T;V) \to \RR$, 
and $\hat{c}:L^2(0,T;V)\times L^2(0,T;V)\to \RR$ together with $\hat{l}_1 : W(0,T)\to\RR$ and $\hat{l}_2 : L^2(0,T;V) \to \RR$ given by
\begin{align*}
    \hat{a}(y,w) &= (y(T),w(T))_H  + \int_0^T(y,w)_V,
    &\hat{b}(w,q) &=  \int_0^T\langle w_t, q \rangle_{V^\star,V}, 
    &\hat{c}(p,q) &= \int_0^T(p,q)_V ,\\
    \hat{l}_1(w) &= (y_0,w(0))_H  + \int_0^T\langle f,w\rangle_{V^\star,V},
    & \hat{l}_2(q) &= \int_0^T \langle f,q\rangle_{V^\star,V} 
\end{align*}
and the problem is given by:\\
Find $(y,p) \in W(0,T)\times L^2(0,T;V)$ such that for all $(w,q) \in W(0,T)\times L^2(0,T;V)$ it holds
\begin{equation}
\tag{\ensuremath{P^s}}
\label{prob:pSaddle}
\begin{aligned}
    \hat{a}(y,w) + \hat{b}(w,p) &= \hat{l}_1(w),\\
    \hat{b}(y,q) - \hat{c}(p,q) &= \hat{l}_2(q).
\end{aligned}    
\end{equation}
The equivalence of \eqref{prob:pSaddle} and \eqref{prob:P} is subject to Lemma~\ref{lem:P-eq-Ps}.

\begin{lemma}
\label{lem:P-eq-Ps}
Let $y \in W(0,T)$ denote the unique solution to \eqref{prob:P}. Then $(y,R(y_t-f)) \in W(0;T)\times L^2(0,T;V)$ denotes a solution to \eqref{prob:pSaddle}. 
Conversely, for the solution $(y,p)$ of \eqref{prob:pSaddle}, $y$ is a solution to \eqref{prob:P}. 
As a consequence, problem \eqref{prob:pSaddle} admits a unique solution.
\end{lemma}
\begin{proof}
Let $y \in W(0,T)$ denote a solution to \eqref{prob:P}. Then it holds
\begin{align}
    \underbrace{(y(T),w(T))_H + \int_0^T(y,w)_V}_{\hat{a}(y,w)} &= 
    \underbrace{(y_0,w(0))_H + \int_0^T \langle f,w\rangle_{V^\star,V}}_{\hat{l}_1(w)}
    + \int_0^T (f,w_t)_{V^\star} - \int_0^T (y_t,w_t)_{V^\star}\\
    &= \hat{l}_1(w) - \int_0^T \langle w_t,R(y_t-f)\rangle_{V^\star,V}.
\end{align}
Setting $p:=R(y_t-f) \in L^2(0,T;V)$ we have
\begin{align}
    \underbrace{\int_0^T(p,q)_V}_{\hat{c}(p,q)} = \underbrace{\int_0^T\langle y_t,q\rangle_{V^\star,V}}_{\hat{b}(y,q)}
    - \underbrace{\int_0^T\langle f,q\rangle_{V^\star,V}}_{\hat{l}_2(q)}\quad \forall q \in L^2(0,T;V),
\end{align}
so that $(y,R(y_t-f))$ satisfies \eqref{prob:pSaddle}. In particular, there exists a solution to 
\eqref{prob:pSaddle}.

On the other hand, let $(y,p) \in W(0,T) \times L^2(0,T;V)$ be a solution to \eqref{prob:pSaddle}. From the second equation of \eqref{prob:pSaddle} we obtain that $p=R(y_t-f)$. Plugging this into the first equation in \eqref{prob:pSaddle} the claim follows. In particular, \eqref{prob:pSaddle} admits a unique solution.
\end{proof}

\begin{remark}
    Existence of solutions to \eqref{prob:pSaddle} can also be proven by showing some LBB-condition. For problem \eqref{prob:pSaddle} in particular the approach from \cite{HongKLP-2022-practicalLBB} is well suited.
\end{remark}


\section{Galerkin approximation.}
\label{GA}
To approximate problem \eqref{prob:P} we propose two approaches. 
In the first approach in Section~\ref{sec:GA-conforming} we investigate a conforming discretization,
while for the second approach, presented in Section \ref{sec:FESP}, we discretize the equivalent saddle point formulation \eqref{prob:pSaddle}.

\subsection{A conforming Galerkin approximation.}
\label{sec:GA-conforming}
For the first approach we use a finite-dimensional subspace $W_d\subset W(0,T)$ to approximate the solution $y$ of \eqref{prob:P} by the unique solution $y_d\in W_d$ of the following conforming variational problem: \\
Find $y_d \in W_d$ such that it holds
\begin{equation}
\tag{\ensuremath{P_d}}
\label{prob:beql_d}
    \begin{aligned}
        b(y_d,w_d) = l(w_d)\quad \forall w_d \in W_d,
    \end{aligned}
\end{equation}
where $b$ and $l$ are defined in \eqref{eq:b_def} and \eqref{eq:l_def}, respectively.

\medskip

For formulation \eqref{prob:beql_d} we have

\begin{theorem}\label{QO}
Problem \eqref{prob:beql_d} admits a unique solution $y_d\in W_d$. Moreover, 
\[\|y-y_d\|_{W(0,T)} \le (1+C_e^2) \inf_{w_d\in W_d}\|y-w_d\|_{W(0,T)},\]
where $C_e$ denotes the constant of the  continuous embedding $W(0,T) \hookrightarrow C([0,T] , H)$.
\end{theorem}

\begin{proof}
Existence and uniqueness of a solution $y_d \in W_d$ follows from the continuity of $l$, and continuity and coercivity of the form $b$. From the definition of $b$ and Galerkin orthogonality it 
with arbitrary $w_d\in W_d$  directly follows that
\[\|y-y_d\|^2_{W(0,T)} \le b(y-y_d,y-y_d) = b(y-y_d,y-w_d) \le  (1+C_e^2) \|y-y_d\|_{W(0,T)}\|y-w_d\|_{W(0,T)}.\]
This directly implies
\[\|y-y_d\|_{W(0,T)} \le (1+C_e^2) \inf_{w_d\in W_d}\|y-w_d\|_{W(0,T)}.\]
\end{proof}
Thus, we have  quasi-optimality of the Galerkin approximation $y_d$ and error estimates can be obtained from e.g.~classical interpolation error estimates, compare e.g. \cite[Thm.~7]{FuehrerKarkulik-SpaceTimeparabolic}.
\medskip

\begin{remark}
Problem~\eqref{prob:beql_d} states the first order necessary and sufficient optimality condition for the minimization problem
\begin{align}
    \min_{w_d \in W_d} \|(w_d)_t + Aw_d - f\|^2_{L^2(0,T;V^{\star})}  + \|y_0 - w_d(0)\|^2_H.
\end{align}
\end{remark}

\subsubsection{Conforming approximation in one spatial dimension.}
\label{sec:1d}
In the formulation \eqref{prob:beql_d} the $V^\star$ norm has to be evaluated. This is possible for certain classes of ansatz functions in one spatial dimension. Here, we for a practical implementation of \eqref{prob:beql_d} use a tensorial space-time finite element  space $W_d$ to approximate $W(0,T)$.
\medskip

Let $I = (0,T]$,  $\Omega \subset \RR^n$,  $V:=H^1_0(\Omega)$ and $H:=L^2(\Omega)$. We use a subdivision  $I_k$ of $I$ with largest interval of length $k$
and a  subdivision $\Omega_h$ of $\Omega$ with largest element of length $h$.
On $I_k$  we define the finite element space 
$K_k = \Span\{\chi_m\mid m=1,\ldots,M\} \subset L^2(I)$ and 
on $\Omega_h$ we define the finite element space 
$V_h = \Span\{\phi_n\mid n=1,\ldots,N\} \subset H^1_0(\Omega)$.
We define the global finite element space  $W_d$ as
\begin{align} 
    W_d = K_k \otimes V_h= 
    \left\lbrace w_d(t,x)=\sum_{m=1}^M\sum_{n=1}^N w_n^m\chi_m(t)\phi_n(x), \;\; w_n^m \in \RR \; \forall n,m \right\rbrace \subset W(0,T).
\end{align}
Thus  we use finite element functions that are separable with respect to time and space.
 We set $d$ as the mesh size parameter that satisfies $d^2 = h^2 + k^2$.
For both $K_k$ and $V_h$ we use piecewise linear and globally continuous finite elements.

\medskip

Assembling the linear algebra representation of \eqref{prob:beql_d} 
requires the evaluation of $L^2(0,T;V^\star)$ inner products of the basis functions of $W_d$.
Using the tensor structure of $W_d$, integrations with respect to time and space are performed independently.
Evaluation with respect to space requires forming $(\phi_n,\phi_m)_{V^\star}$ for all test functions $\phi_n,\phi_m$, $n,m=1,\ldots,N$.
These integrals are evaluated using the corresponding analytical Riesz representer $R\phi_n$ of $\phi_n$ via the identity $(\phi_n,\phi_m)_{V^\star} = (R\phi_n,\phi_m)_H$. This leads to dense matrices of size $N\times N$. The resulting linear algebra representation of \eqref{prob:beql_d} 
is a tridiagonal block matrix of size $M N \times M N$, 
that contains dense blocks of size $N \times N$ on the three main diagonals. 
For the purpose of this study, these systems are solved using direct solvers.

How to compute $R\phi_n$ is shown in the next Section~\ref{sec:Riesz1D}. 
Numerical results are shown in Section~\ref{sec:num-1D}.

\subsubsection{Riesz representer of piecewise linear, continuous functions.}
\label{sec:Riesz1D}
The numerical scheme \eqref{prob:beql_d} requires the evaluation of the exact Riesz representer of basis functions $\phi_n \in V_h$. In the following we show how this representer for piecewise linear and globally continuous basis functions in one space dimension can be evaluated analytically. To begin with, let $\Omega=(a,b) \subset \RR$. 
We consider a grid with $N$ vertices $a=x_1 < \cdots < x_N=b$. 
Let $\phi_n$ denote the $n$-th piecewise linear and globally continuous basis function, defined by $\phi_n(x_j) = \delta_{nj}$, for $n,j = 1,\ldots,N$.
Since $w_n:= R \phi_n$ solves $-w_n''=\phi_n$ in $(a,b)$, $w_n(a)=w_n(b)=0$, a short calculation shows 
\begin{equation}
    R \phi_1 (x) = 
    \begin{cases} 
    \beta_1 (x - x_N) - \frac{(x-x_{2})^3}{6(x_1 - x_{2})}, & x\in \left[x_{1}, x_{2}\right) \\
    \beta_1 (x - x_N), & x\in \left[x_2, x_{N}\right]
    \end{cases}, 
    \qquad \beta_1 := \frac{(x_1 - x_2)^2}{6(x_1-x_N)}
\end{equation}
and
\begin{equation}
    R \phi_N (x) = 
    \begin{cases} 
    \alpha_N (x - x_1) , & x\in \left[x_{1}, x_{N-1}\right) \\
    \alpha_N (x - x_1) - \frac{(x-x_{N-1})^3}{6(x_N - x_{N-1})}, & x\in \left[x_{N-1}, x_{N}\right]
    \end{cases}, 
    \qquad \alpha_N := \frac{(x_N - x_{N-1})^2}{6(x_N-x_1)}.
\end{equation}
and for  $2\leq n\leq N-1$
\begin{align}
    R \phi_n (x) = \begin{cases} \alpha_n (x - x_1), & x\in \left[x_1, x_{n-1}\right),\\
    \alpha_n (x - x_1) - \frac{(x-x_{n-1})^3}{6(x_n - x_{n-1})}, & x\in \left[x_{n-1}, x_{n}\right), \\
    \beta_n (x - x_N) - \frac{(x-x_{n+1})^3}{6(x_n - x_{n+1})}, & x\in \left[x_{n}, x_{n+1}\right), \\
    \beta_n (x - x_N), & x\in \left[x_{n+1}, x_{N}\right],
    \end{cases}
\end{align}
where 
\begin{align}
    \alpha_n &:= \frac{x_{n-1}-x_{n+1}}{x_N - x_1} \left(\frac{1}{6}(x_{n+1} - 2x_n + x_{n-1}) 
    + \frac{1}{2} (x_{n} - x_{N}) \right), \\
    \beta_n &:= \frac{x_{n-1}-x_{n+1}}{x_N - x_1} \left(\frac{1}{6}(x_{n+1} - 2x_n + x_{n-1}) 
    + \frac{1}{2} (x_{n} - x_{1}) \right) 
    = \alpha_n + \frac{1}{2}(x_{n-1} - x_{n+1}).
\end{align}
In Figure~\ref{fig:FEdisc:Riesz} we show a basis function together with its Riesz representer.

\begin{figure}
    \centering

\begin{tikzpicture}
\pgfmathsetmacro{\xo}{0};
\pgfmathsetmacro{\xN}{9};

\begin{axis}
[
width=0.9\textwidth
,height=3cm
,xmin=\xo
,xmax=\xN
,ymin=0
,ymax=2.0
,xtick={2,3,4}
,xticklabels={$x_{i-1}$,$x_i$,$x_{i+1}$}
]



\pgfmathsetmacro{\xm}{2.0}
\pgfmathsetmacro{\xi}{3.0}
\pgfmathsetmacro{\xp}{4.0}
\pgfmathsetmacro{\alpha}{(\xm-\xp)/(\xN-\xo)*((\xi-\xN)/2}
\pgfmathsetmacro{\beta}{\alpha+(\xm-\xp)/2}
\addplot[color=black,dotted,thick] table[row sep=\\]
{
2.0 0.0\\
3.0 1.0\\
4.0 0.0\\
};
\addlegendentry{$\phi_i$}

\addplot[samples=2,domain=0:2]{\alpha*(x-\xo)};
\addplot[samples=12,domain=2:3]{\alpha*(x-\xo)-((x-\xm)^3)/6};
\addplot[samples=12,domain=3:4]{\beta*(x-\xN)+((x-\xp)^3)/6};
\addplot[samples=2,domain=4:9]{\beta*(x-\xN)};
\addlegendentry{$R\phi_i$}



\end{axis}

\end{tikzpicture}
    \caption{The piecewise linear basis function $\phi_n$ together with its
    Riesz representers $R\phi_n$.}
    \label{fig:FEdisc:Riesz}
\end{figure}

\begin{remark}
The exact calculation of $R\phi_n$ in higher-dimensional cases can be performed with the help of Green's function. For a rectangle we refer to \cite[§~9.7]{Myint} together with \cite[§~11.4(1)]{Myint} and for the unit disk to \cite[§~11.5]{Myint}. 
For more general domains alternatives to the exact evaluation of $R$ could be considered. This is subject to future work.
\end{remark}


\subsection{A finite element approximation in general spatial dimension.}
\label{sec:FESP}
For a fully practical approach, also the Riesz lift $R$ must be discretized. In the following we introduce a fully discrete approximation of \eqref{prob:P} using a discrete approximation of $R$.
We investigate this formulation with the help of a discrete version of the saddle point problem \eqref{prob:pSaddle}.
 To begin with we in the following let $I = (0,T]$,  $\Omega \subset \RR^{n}$, $n\in \{1,2,3\}$, $V:=H^1_0(\Omega)$ and $H:=L^2(\Omega)$.

\medskip

We recall that we understand the Riesz lift $R:V^\star\to V$ as operator $R:L^2(0,T;V^\star) \to L^2(0,T;V)$ by acting identical with respect to time and define its discrete approximation 
$R_d : L^2(0,T;V^\star) \to Q_d$ for arbitrary  $\phi \in L^2(0,T;V^{\star})$ by the equation
\begin{align}
\label{eq:SP_def_Riesz_discr}
    \int_0^T(R_d\phi,q_d)_V = \int_0^T \langle \phi,q_d\rangle_{V^\star,V} \quad \forall q_d \in Q_d.
\end{align}
With this definition 
we define the non conforming discrete approximation of \eqref{prob:P} as:\\
Find $y_d \in W_d$ such that for all $w_d \in W_d$ it holds
\begin{align}
\tag{\ensuremath{P_{dd}}}
\label{eq:bd-eq-ld}
    b_d(y_d,w_d) = l_d(w_d),
\end{align}
where
\begin{align}
    b_d(y_d,w_d) &:=  (y_d(T),w_d(T))_H + \int_0^T(y_d,w_d)_V + \int_0^T(R_d(y_d)_t,R_d(w_d)_t)_V,
    \label{eq:bd}\\
    l_d(w_d) &:= (y_0,w_d(0))_H + \int_0^T\langle f,w_d\rangle_{V^\star,V} + \int_0^T \langle f,R_d(w_d)_t\rangle_{V^\star,V} .
\end{align}

\begin{remark}
    Let $\phi \in L^2(0,T;V^\star)$ be arbitrary.
\begin{itemize}
    \item Since $\int_0^T(R\phi,q_d)_V = \int_0^T\langle \phi,q_d\rangle_{V^\star,V}$, Galerkin orthogonality $\int_0^T (R\phi - R_d \phi,q_d)_V = 0$ is satisfied from which we directly obtain the stability estimate
 \begin{align*}
\|R_d \phi\|_{L^2(0,T;V)} \leq \|R\phi\|_{L^2(0,T;V)} \leq \|\phi\|_{L^2(0,T;V^\star)}.     
 \end{align*}
 \item Moreover, from the identity  $\int_0^T(R\phi-R_d\phi,R_d\phi-R\phi+R\phi)_{V}=0$ we obtain 
 \begin{align*}
     \|R\phi-R_d\phi\|_{L^2(0,T;V)} \leq\|R\phi\|_{L^2(0,T;V)} = \|\phi\|_{L^2(0,T;V^\star)}.
 \end{align*}

\end{itemize}

\end{remark}

\begin{lemma}
\label{lem:bdld_exsol}
    The variational problem
 \eqref{eq:bd-eq-ld}
admits a unique solution $y_d\in W_d$ for any  fixed $d>0$.
\end{lemma}
\begin{proof}
    Since $W_d$ is finite dimensional, all norms on $W_d$ are equivalent. Let $W_d$ be equipped with the norm induced by $W(0,T)$. Then $b_d$ is bilinear, continuous and coercive on $W_d$, and $l_d$ is continuous on $W_d$. Thus, problem \eqref{eq:bd-eq-ld} admits a unique solution $y_d\in W_d$.
\end{proof}

Since we do not gain bounds from Lemma~\ref{lem:bdld_exsol} that are uniform in $d$, we introduce an equivalent formulation using the saddle point formulation \eqref{prob:pSaddle}.

\medskip

We introduce finite element spaces $W_d \subset W(0,T)$ and $Q_d \subset L^2(0,T;V)$ that provide sufficient approximation properties for $d \to 0$.

\begin{assumption}$ $
\begin{enumerate}
    \item We assume that $W_d$ and $Q_d$ are chosen such that the compatibility condition
\begin{align}
\label{compC}
    (w_d)_t \in Q_d \text{ for all } w_d \in W_d
\end{align}
holds.     
\item 
The $L^2$-projection $\Pi_d : L^2(0,T;L^2(\Omega)) \to Q_d$ is defined by
\begin{align}
    \int_0^T( \Pi_d v,q_d)_H = \int_0^T(v,q_d)_H\quad \text{for all } q_d \in Q_d.
\end{align}
We assume, that the finite element mesh with respect to space (resp. space-time) is such, that the $L^2$-projection is $L^2(0,T;V)$-stable,  
i.e., there is $C>0$, independent of $d$, such that for all $v \in L^2(0,T;V)$
\begin{align}
 \| \Pi_d v\|_{L^2(0,T;V)} \leq C \|v\|_{L^2(0,T;V)}. 
\end{align}
This is satisfied for example on quasi-uniform meshes or shape regular meshes under suitable additional conditions \cite{H1L2stability-Bramble2001,H1L2stability-Bank2013}.
\end{enumerate}

\end{assumption}

We define the finite element approximation of \eqref{prob:pSaddle} as:\\
Find $(y_d,p_d) \in W_d \times Q_d$
such that 
\begin{equation}
\tag{\ensuremath{P^s_d}}
\label{prob:pSaddle_disc}
\begin{aligned}
    \hat{a}(y_d,w_d) + \hat{b}(w_d,p_d) &= \hat{l}_1(w_d) \quad \forall w_d \in W_d,\\
    \hat{b}(y_d,q_d) - \hat{c}(p_d,q_d) &= \hat{l}_2(q_d)\quad \forall q_d \in  Q_d.
\end{aligned}
\end{equation}

\begin{lemma} 
\label{lem:bdld-eq_SPd}
Let $(y_d,p_d)$ denote a solution to \eqref{prob:pSaddle_disc}. 
Then $y_d$ is a solution to \eqref{eq:bd-eq-ld}.
Vice versa, if $y_d$ denotes the solution to \eqref{eq:bd-eq-ld}, then $(y_d,R_d( (y_d)_t - f))\in W_d\times Q_d$ is the solution to \eqref{prob:pSaddle_disc}. 
\end{lemma}

\begin{proof}
Let $(y_d,p_d) \in W_d \times Q_d$ denote a solution to \eqref{prob:pSaddle_disc}. 
Then especially we have $\hat{c}(p_d,q_d) = \hat{b}(y_d,q_d)-\hat{l}_2(q_d)$ and thus by construction and by \eqref{eq:SP_def_Riesz_discr} $p_d = R_d ( (y_d)_t - f) \in Q_d$. 
Inserting this into $\hat{b}(w_d,p_d)$ we obtain
\begin{align}
   \hat{b}(w_d, R_d( (y_d)_t - f)) &=  \int_0^T \langle (w_d)_t, R_d( (y_d)_t - f)\rangle_{V^\star,V} \\
   & = \int_0^T \langle (w_d)_t, R_d(y_d)_t\rangle_{V^\star,V}
   - \int_0^T \langle (w_d)_t, R_d f\rangle_{V^\star,V}\\
   & = \int_0^T ( R_d (w_d)_t,  R_d (y_d)_t )_{V} -\int_0^T ( R_d (w_d)_t, R_d f )_{V}\\
    &= \int_0^T (R_d (y_d)_t,  R_d(w_d)_t)_V
    -\int_0^T \langle f,R_d (w_d)_t \rangle_{V^\star,V}.
\end{align}
In summary $y_d$ satisfies
\[b_d(y_d,w_d)=l_d(w_d) \text{ for all } w_d\in W_d.\]

Finally, if $y_d\in W_d$ is the unique solution to \eqref{eq:bd-eq-ld}, then $(y_d,p_d)$ with $p_d:= R_d ( (y_d)_t - f) \in Q_d$ solves \eqref{prob:pSaddle_disc} and is its unique solution.

\end{proof}

From Lemma~\ref{lem:bdld-eq_SPd} we obtain, that for any fixed $d>0$ problem \eqref{prob:pSaddle_disc} admits a unique solution which also is the solution to \eqref{eq:bd-eq-ld}. In the following we consider the limit $d\to 0$ and derive uniform bounds on the solution to \eqref{prob:pSaddle_disc}.

\begin{lemma}
\label{lem:ydpd-bndd}
Let $(y_d,p_d) \in W_d \times Q_d$ denote the solution to \eqref{prob:pSaddle_disc}.
Then
\begin{align}
   \frac{1}{2}\|y_d(T)\|_H^2 + \frac{3}{8}\|y_d\|^2_{L^2(0,T;V)} + \frac{3}{8}\|p_d\|^2_{L^2(0,T;V)} \leq \frac{1}{2}\|y_0\|_H^2 + 4\|f\|_{L^2(0,T;V^\star)}^2.
\end{align}
\end{lemma}
\begin{proof}
Using $q=p_d$ and $w = y_d$ and subtracting the equations we observe
\begin{align}
    \|y_d(T)\|_H^2 &+ \int_0^T (y_d,y_d)_V 
    + \int_0^T \langle(y_d)_t,p_d\rangle_{V^\star,V}
    - \int_0^T \langle(y_d)_t,p_d\rangle_{V^\star,V}
    +\int_0^T (p_d,p_d)_V\\
    &=(y_0,y_d(0))_H+\int_0^T \langle f,y_d\rangle_{V^\star,V}
     - \int_0^T \langle f,p_d\rangle_{V^\star,V}
\end{align}
and thus 
\begin{align}
   \|y_d(T)\|_H^2 +  \|y_d\|^2_{L^2(0,T;V)} + \|p_d\|_{L^2(0,T;V)}^2 
   =    (y_0,y_d(0))_H
   + \int_0^T \langle f,y_d\rangle_{V^\star,V} 
   - \int_0^T \langle f,p_d\rangle_{V^\star,V}. \label{eq:ydpdbnd-1}
\end{align}
To estimate $(y_0,y_d(0))_H$ we proceed
\begin{align}
\frac{1}{2} \|y_d(T)\|_H^2 - \frac{1}{2}\|y_d(0)\|_H^2 
&= \frac{1}{2}\int_0^T \frac{d}{dt}\|y_d\|_H^2 
= \int_0^T \langle (y_d)_t,y_d\rangle_{V^\star,V} \\
&= \int_0^T ((y_d)_t,y_d)_H = \int_0^T((y_d)_t,\Pi_d y_d)_H
=\int_0^T\langle (y_d)_t, \Pi_d y_d\rangle_{V^\star,V}\\
&=\int_0^T (p_d, \Pi_d y_d)_V 
+  \int_0^T\langle f, \Pi_d y_d\rangle_{V^\star,V}.
\end{align}
Here we used that $(y_d)_t \in Q_d \subset L^2(0,T;L^2(\Omega))$ to substitute 
$y_d \in W_d$
by $\Pi_d y_d \in Q_d$. Since $p_d \in Q_d$ the projection in the first term can be skipped and we thus obtain
\begin{align}
    \frac{1}{2} \|y_d(0)\|_H^2 
    &= \frac{1}{2}\|y_d(T)\|_H^2 - \int_0^T(p_d,y_d)_V 
    -  \int_0^T\langle f, \Pi_d y_d\rangle_{V^\star,V} . \label{eq:ydpdbnd-2}
\end{align}
Proceeding from \eqref{eq:ydpdbnd-1} and inserting \eqref{eq:ydpdbnd-2} we get 
\begin{align}
\begin{split}
  \|y_d(T)\|_H^2 +  \|y_d\|^2_{L^2(0,T;V)} 
  &+ \|p_d\|_{L^2(0,T;V)}^2 \\
  &=
   (y_0,y_d(0))_H + 
    \int_0^T \langle f,y_d\rangle_{V^\star,V} - \int_0^T \langle f,p_d\rangle_{V^\star,V}\end{split}\\
    & \leq \frac{1}{2}\|y_0\|_H^2 + \frac{1}{2}\|y_d(0)\|_H^2 
     + \int_0^T \langle f,y_d\rangle_{V^\star,V} - \int_0^T \langle f,p_d\rangle_{V^\star,V}\\
     \begin{split}
     &= \frac{1}{2}\|y_0\|_H^2  + \frac{1}{2}\|y_d(T)\|_H^2
     - \int_0^T(p_d,y_d)_V -  \int_0^T\langle f, \Pi_d y_d\rangle_{V^\star,V} \\
     &\hphantom{=}+ \int_0^T \langle f,y_d\rangle_{V^\star,V} - \int_0^T \langle f,p_d\rangle_{V^\star,V} \end{split}\\
     \begin{split}
     &\leq \frac{1}{2}\|y_0\|_H^2  + \frac{1}{2}\|y_d(T)\|_H^2
     + \frac{1}{2}\|p_d\|^2_{L^2(0,T;V)} + \frac{1}{2}\|y_d\|^2_{L^2(0,T;V)}\\
     &\hphantom{=}+ \int_0^T \langle f,y_d\rangle_{V^\star,V} - \int_0^T \langle f,p_d\rangle_{V^\star,V}
     -  \int_0^T (R_d f, \Pi_d y_d)_V. 
     \end{split}
\end{align}
In the last step we use the definition of $R_d$ together with $\Pi_d y_d \in Q_d$. 
Since $R_df \in Q_d$ we can omit $\Pi_d$ on the last term.
We conclude
\begin{align}
     \frac{1}{2}\|y_d(T)\|_H^2
     &+  \frac{1}{2}\|y_d\|^2_{L^2(0,T;V)} 
  + \frac{1}{2}\|p_d\|_{L^2(0,T;V)}^2 \\
    &\leq \frac{1}{2}\|y_0\|_H^2  
    + \int_0^T \langle f,y_d\rangle_{V^\star,V} - \int_0^T (R_d f, y_d)_V
    - \int_0^T \langle f,p_d\rangle_{V^\star,V} \\
    & \leq \frac{1}{2}\|y_0\|_H^2 
    +\|Rf-R_df\|_{L^2(0,T;V)} \|y_d\|_{L^2(0,T;V)}
     + \|f\|_{L^2(0,T;V^\star)}\|p_d\|_{L^2(0,T;V)}, \\
    & \leq \frac{1}{2}\|y_0\|_H^2 
    +\|f\|_{L^2(0,T;V^\star)}\|y_d\|_{L^2(0,T;V)}
     + \|f\|_{L^2(0,T;V^\star)}\|p_d\|_{L^2(0,T;V)},\\
    & = \frac{1}{2}\|y_0\|_H^2 
    +2 \|f\|_{L^2(0,T;V^\star)} \frac{1}{2} \|y_d\|_{L^2(0,T;V)}
     + 2\|f\|_{L^2(0,T;V^\star)} \frac{1}{2}\|p_d\|_{L^2(0,T;V)},
\end{align}
where we use the stability of the discrete Riesz lift.
By using Young's inequality we finish with
\begin{align}
    \frac{1}{2}\|y_d(T)\|_H^2
     + \frac{3}{8}\|y_d\|^2_{L^2(0,T;V)} 
  + \frac{3}{8}\|p_d\|_{L^2(0,T;V)}^2 
     \leq \frac{1}{2}\|y_0\|_H^2 
     + 4\|f\|_{L^2(0,T;V^\star)}^2.
\end{align}
\end{proof}

To show uniform boundedness of the solution $(y_d,p_d)$ to \eqref{prob:pSaddle_disc} in $W(0,T)\times L^2(0,T;V)$ a bound on $(y_d)_t$ is missing. 

\begin{lemma}
\label{lem:ydt-bndd}
Let the $L^2$-projection $\Pi_d : L^2(0,T;V) \to Q_d$ be stable with respect to the $L^2(0,T;V)$-norm. Then there exists $C>0$ independent of $d$ such that
\begin{align}
    \|(y_d)_t\|_{L^2(0,T;V^\star)} 
    \leq C \left( \|p_d\|_{L^2(0,T;V)} + \|f\|_{L^2(0,T;V^\star)} \right).
\end{align}
\end{lemma}
\begin{proof}
We use $q_d = \Pi_d(R (y_d)_t) \in Q_d $ as test function in \eqref{prob:pSaddle_disc} and obtain 
\begin{align}
    \hat{b}(y_d,\Pi_d(R(y_d)_t) ) 
    &= \int_0^T\langle f, \Pi_d(R(y_d)_t )\rangle_{V^\star,V} 
    +  \hat{c}(p_d, \Pi_d(R(y_d)_t) ) .
\end{align}
For the left hand side we write using $(y_d)_t \in Q_d$ and the properties of the
$L^2$-projection
\begin{align}
    \hat{b}(y_d, \Pi_d (R(y_d)_t)) 
    &= \int_0^T\langle (y_d)_t, \Pi_d (R(y_d)_t)\rangle_{V^\star,V} 
    = \int_0^T( (y_d)_t, \Pi_d(R(y_d)_t))_H \\
&= \int_0^T( (y_d)_t,(R(y_d)_t))_H
= \int_0^T  \langle (y_d)_t,R(y_d)_t\rangle_{V^\star,V}
= \| (y_d)_t\|^2_{L^2(0,T;V^\star)} .
\end{align}
On the other hand we have the existence of $C>0$ independent of $d$ 
using the $L^2(0,T;V)$-stability of $\Pi_d$  such that
\begin{align}
    \hat{c}(p_d, \Pi_d(R(y_d)_t)) 
    = \int_0^T (p_d, \Pi_d (R(y_d)_t))_V
    & \leq \int_0^T \|p_d\|_V\|\Pi_d (R(y_d)_t)\|_V\\
    &\leq C \int_0^T \|p_d\|_V \|R(y_d)_t\|_V 
     \leq C \|p_d\|_{L^2(0,T;V)} \|(y_d)_t\|_{L^2(0,T;V^\star)}
\end{align}
and using the same arguments
\begin{align}
    \int_0^T \langle f, \Pi_d (R(y_d)_t) \rangle_{V^\star,V}
    &\leq C\|f\|_{L^2(0,T;V^\star)}  \|(y_d)_t\|_{L^2(0,T;V^\star)} .
\end{align}
Combining the estimates we obtain
\begin{align}
\|(y_d)_t\|_{L^2(0,T;V^\star)} \leq C \left(\|p_d\|_{L^2(0,T;V)} + \|f\|_{L^2(0,T;V^\star)}\right).
\end{align}

\end{proof}

\begin{corollary}
By combining the results of Lemma~\ref{lem:ydpd-bndd} and Lemma~\ref{lem:ydt-bndd}
there exists $C>0$ independent of $d$ such that
\begin{align}
   \|y_d(T)\|_H +   \|y_d\|_{W(0,T)} + \|p_d\|_{L^2(0,T;V)} \leq C.
\end{align}
\end{corollary}

\begin{theorem}
The sequences $(y_d)_d$ and $(p_d)_d$ weakly converge in $W(0,T)\times L^2(0,T;V)$ to the solution $(y,p)$ of problem \eqref{prob:pSaddle}.
\end{theorem}
\begin{proof}
Since $(y_d)_d$ and $(p_d)_d$ are bounded, there exist weakly convergent subsequences $(y_{d_l})_l, (p_{d_l})_l$ with the property
\begin{align}\label{eq:proof:weaksubseq}
    (y_{d_l},p_{d_l}) \rightharpoonup (y,p) \quad \text{for} \quad d_l\rightarrow 0.
\end{align}
We show that $(y,p)$ from \eqref{eq:proof:weaksubseq} satisfies \eqref{prob:pSaddle} for all $(w,q) \in W(0,T)\times L^2(0,T;V)$. To begin let $(w,q)$ be given together with sequences $(w_{d_l})_l$ and $(q_{d_l})_l$ converging strongly to $(w,q)$ in $W(0,T)\times L^2(0,T;V)$ for $l\rightarrow \infty$. By \eqref{eq:proof:weaksubseq} we have that
\begin{align}
    \hat{a}(y_{d_l},w) + \hat{b}(w,p_{d_l}) + \hat{b}(y_{d_l},q) - \hat{c}(p_{d_l},q) &\longrightarrow \hat{a}(y,w) + \hat{b}(w,p) + \hat{b}(y,q) - \hat{c}(p,q) \qquad (d_l \rightarrow 0).
\end{align}
On the other hand for the left hand side we observe that
\begin{align}
    \begin{split}
        \hat{a}(y_{d_l},w) + \hat{b}(w,p_{d_l}) + \hat{b}(y_{d_l},q) - \hat{c}(p_{d_l},q) &= \hat{a}(y_{d_l},w_{d_l}) + \hat{b}(y_{d_l},q_{d_l}) \\&\hphantom{=}+  \hat{b}(w_{d_l},p_{d_l}) - \hat{c}(p_{d_l},q_{d_l})     \end{split}\tag{$\dagger$}\label{eq:proof:up}\\\begin{split}&\hphantom{=}+ \hat{a}(y_{d_l},w-w_{d_l}) + \hat{b}(y_{d_l},q-q_{d_l}) \\&\hphantom{=}+  \hat{b}(w-w_{d_l},p_{d_l}) - \hat{c}(p_{d_l},q-q_{d_l}) .\end{split}\tag{$\ddagger$}\label{eq:proof:down} 
\end{align}
For the part \eqref{eq:proof:up} we obtain
\begin{align}
    \hat{a}(y_{d_l},w_{d_l}) + \hat{b}(y_{d_l},q_{d_l}) +  \hat{b}(w_{d_l},p_{d_l}) - \hat{c}(p_{d_l},q_{d_l}) &= \hat{l}_1(w_{d_l}) + \hat{l}_2(q_{d_l}) \longrightarrow \hat{l}_1(w) + \hat{l}_2(q) \qquad (d_l \rightarrow 0).
\end{align}
Moreover for \eqref{eq:proof:down} we have
\begin{align}
    \hat{a}(y_{d_l},w-w_{d_l}) + \hat{b}(y_{d_l},q-q_{d_l}) +  \hat{b}(w-w_{d_l},p_{d_l}) - \hat{c}(p_{d_l},q-q_{d_l}) \longrightarrow 0 \qquad (d_l \rightarrow 0).
\end{align}
Thus, for $d_l \rightarrow 0$ the limit $(y,p)$ solves \eqref{prob:pSaddle}. Since all other weakly convergent subsequences would converge to the same limits we have $y_d\rightharpoonup y$ and $p_d\rightharpoonup p$ for $d\rightarrow 0$.

\end{proof}

\begin{remark}
    Let us note that we do not impose a LBB condition for the solution of the saddle point problem \eqref{prob:pSaddle_disc}. Instead, we require the compatibility condition \eqref{compC} for the discrete spaces $W_d$ and $Q_d$ (which of course might be interpreted as a LBB condition).
\end{remark}

\subsection{Tensorial finite element discretization of the saddle point problem.}\label{subsec:FEDiscr}
Our formulation and the previous results hold true in both cases, where a separate discretization of space and time are used or a full space-time discretization. The latter is meaningful especially in the context of shape optimization with time-dependent shapes. Results on that will be published elsewhere. For the sake of simplicity, we here focus on tensorial space-time finite elements to approximate $W(0,T)$ and $L^2(0,T;V)$.

We use a grid $I_k$ in time
and a conforming triangulation $\Omega_h$ in space.
We set $d$ as the mesh size parameter that satisfies $d^2 = h^2 + k^2$, where $k$ denotes the mesh size of $I_k$ and $h$ denotes the mesh size of $\Omega_h$.  
On $I_k$ and $\Omega_h$ we define finite element spaces 
\begin{align}
J_k &= \Span\{\psi_p\mid p=1,\ldots,P\} \subset L^2(I),\\
K_k &= \Span\{\chi_m\mid m=1,\ldots, M, (\chi_m)_t \in J_k\} \subset L^2(I),\\
V_h &= \Span\{\phi_n\mid n=1,\ldots,N\} \subset V
\end{align}
and define global finite element spaces  $Q_d$ and $W_d$  as
\begin{align} 
    Q_d = J_k \otimes V_h &=
    \left\lbrace q_d(t,x)=\sum_{p=1}^P\sum_{n=1}^N q_n^p\psi_p(t)\phi_n(x), \;\; q_n^p \in \RR \; \forall n,p \right\rbrace \subset L^2(0,T;V) \label{eq:Qd},\\
    W_d = K_k \otimes V_h &= 
    \left\lbrace w_d(t,x)=\sum_{m=1}^M\sum_{n=1}^N w_n^m\chi_m(t)\phi_n(x), \;\; w_n^m \in \RR \; \forall n,m \right\rbrace \subset W(0,T). \label{eq:Wd}
\end{align}
For both $K_k$ and $V_h$ we use piecewise linear and globally continuous finite elements and for $J_k$ we use piecewise constant and discontinuous finite elements.
Thus especially we use finite element functions that are separable with respect to time and space which will be exploited in generating and solving the resulting linear problems.

To state the final linear system, let $M_t, T_t \in \RR^{M \times M}$, $M_t^q \in \RR^{P\times P}$, $Z_t \in \RR^{P \times M}$ and $M_x, A_x \in \RR^{N \times N}$ be given by
\begin{alignat}{6}
& (M_t)_{ij} &&= \int_0^T \chi_j \chi_i, \qquad && (T_t)_{ij} &&= \chi_j(T)\chi_i(T), \qquad && (M_t^q)_{ij} &&= \int_0^T \psi_j \psi_i, \\
& (Z_t)_{ij} &&= \int_0^T (\chi_j)_t \psi_i, \qquad && (M_x)_{ij} &&= (\phi_j, \phi_i)_H, \qquad && (A_x)_{ij} &&= (\phi_j, \phi_i)_V .
\end{alignat}
Note that the mass $M_t^q$ is different to the mass $M_t$. Moreover we introduce $F_1^t, R_0^t \in \RR^M$, $F_2^t \in \RR^P$ and $F^x, R_0^x \in \RR^N$ by
\begin{alignat}{6}
& (F_1^t)_{i} &&= \int_0^T f^t \chi_i, \qquad && (R_0^t)_{i} &&= \chi_i(0), \qquad && (F_2^t)_{i} &&= \int_0^T f^t \psi_i, \\
& (F^x)_{i} &&= \langle f^x, \phi_i \rangle_{V^\star, V}, \qquad && (R_0^x)_{i} &&= (y_0, \phi_i)_H . \qquad &&
\end{alignat}
Let $y\in\RR^{PN},p\in\RR^{MN}$ denote the degrees of freedom of $y_d$ and $p_d$.
On the linear algebra level \eqref{prob:pSaddle_disc} then reads as a linear system of the form
\begin{align}
    \begin{pmatrix}
    A & B^T\\ B & -C
    \end{pmatrix}
    \begin{pmatrix}
    y \\
    p
    \end{pmatrix}
    = 
    \begin{pmatrix}
    F_1 \\ F_2
    \end{pmatrix} ,
\end{align}
where
\begin{alignat}{6}
& A &&= T_t \otimes M_x + M_t \otimes A_x, \qquad && B &&=  Z_t \otimes M_x, \qquad && C &&= M_t^q \otimes A_x, \\
& F_1 &&= F_1^t \otimes F^x +  R_0^t \otimes R_0^x, \qquad && && && F_2 &&= F_2^t \otimes F^x,
\end{alignat}
Here $\otimes$ stands for the Kronecker product of two matrices.
For details we refer to \cite{Golub}.

\subsection{Derivation of a preconditioner.}
\label{sec:prec}
For a numerical treatment we reorder the system and obtain
\begin{align}
    \begin{pmatrix}
    C & -B\\ -B^T & -A
    \end{pmatrix}
    \begin{pmatrix}
    p \\
    y
    \end{pmatrix}
    = 
    \begin{pmatrix}
    -F_2 \\ -F_1
    \end{pmatrix} .
\end{align}
An ideal left sided preconditioner now is given by \cite{BenziGolubLiesen}
\begin{align}
    P = 
    \begin{pmatrix}
    C & -B\\
    0 & -S
    \end{pmatrix},
    \quad S = A+B^TC^{-1}B ,
    \label{eq:P_left}
\end{align}
where $S$ is the Schur complement. 
For a practical preconditioner we use approximations of $C$ and $S$, which we specify in the following. 

If $J_k$ is spanned by piecewise constant finite elements, $M_t^q$ is a diagonal matrix and $C = M_t^q \otimes A_x$ is blockdiagonal, with $A_x$ on the diagonal. Since $A_x$ represents a diffusion operator, preconditioning using multigrid is suitable.
For other finite elements, we approximate $M_t^q$ by its diagonal and proceed as before.

For the approximation of $S$ we observe that
\begin{align}
    S &= T_t \otimes M_x + M_t \otimes A_x + (Z_t^T \otimes M_x)(M_t^q \otimes A_x)^{-1}(Z_t\otimes M_x)\\
    &= T_t \otimes M_x + M_t \otimes A_x + Z_t^T(M_t^q)^{-1}Z_t \otimes M_x A_x^{-1}M_x . \label{eq:matrixSchur}
\end{align}
Let $\Lambda_t,Q_t \in \RR^{M\times M}$ be obtained from the generalized eigenvalue problem $Z_t^T (M_t^q)^{-1} Z_t \Vec{x}_i = \lambda_i M_t \Vec{x}_i$, i.e.~$\Lambda_t$ is a diagonal matrix holding the generalized eigenvalues on its diagonal, while $Q_t$ contains the corresponding generalized eigenvectors as columns. These matrices satisfy
\begin{align}
    Q_t^T (Z_t^T (M_t^q)^{-1} Z_t) Q_t = \Lambda_t \quad \text{and} \quad Q_t^T M_t Q_t = I_t ,
\end{align}
where  $I_t \in \RR^{M}$ is the identity matrix on $\RR^{M}$. 

In a first approximation step we introduce $D_t = \operatorname{diag}(Q_t^T T_t Q_t) =: \operatorname{diag}(d_1,d_2,\ldots,d_M)$ and 
$I_x \in \RR^{N\times N}$ as identity matrix on $\RR^{N}$ and proceed 
\begin{align}
    S &= (Q_t^{-T} \otimes I_x)(Q_t^T T_t Q_t \otimes M_x + I_t \otimes A_x + \Lambda_t \otimes M_x A_x^{-1}M_x)(Q_t^{-1} \otimes I_x)\\
    &\approx(Q_t^{-T} \otimes I_x)(D_t \otimes M_x + I_t \otimes A_x + \Lambda_t \otimes M_x A_x^{-1}M_x)(Q_t^{-1} \otimes I_x).
\end{align}
Now the block
$D_t \otimes M_x + I_t \otimes A_x + \Lambda_t \otimes M_x A_x^{-1}M_x$ is blockdiagonal and we prepare a preconditioner for each of the diagonal blocks individual. The $i$-th diagonal block is given by 
$d_i M_x + A_x + \lambda_i M_xA_x^{-1}M_x$ and we introduce the approximation
\begin{align}
    d_i M_x + A_x + \lambda_i M_xA_x^{-1}M_x 
    \approx A_x + 2\sqrt{\lambda_i} M_x +  \lambda_i M_xA_x^{-1}M_x
    =(A_x + \sqrt{\lambda_i}M_x)A_x^{-1}(A_x + \sqrt{\lambda_i}M_x).
\end{align}
following \cite{PearW-2012-PrecondSchur}.
For solving these systems, again multigrid methods are used \cite{Xu_Zikatanov_2017}.

\begin{remark}
\label{rem:stiffnessTime}
In the standard setting, where piecewise linear and globally continuous finite elements (CG 1) are spanning $K_k$ and piecewise constant elements (DG 0) are used for $J_k$, it holds that
\begin{align}
   (Z_t^T (M_t^q)^{-1} Z_t)_{ij} = \int_0^T (\chi_j)_t (\chi_i)_t .
\end{align}
Moreover 
the Schur complement \eqref{eq:matrixSchur} can be considered as the discrete representation of \eqref{eq:bd}.
\end{remark}

\begin{remark}
Since we solve a symmetric saddle point problem, using a symmetric solver like MINRES is an option. This requires a symmetric positive definite preconditioner. In this situation instead of \eqref{eq:P_left} we
consider
\begin{align}
    P_\text{sym} = 
    \begin{pmatrix}
    C & 0\\
    0 & S
    \end{pmatrix}, 
    \label{eq:P_sym}
\end{align}
as preconditioner and proceed with the same approximation steps as before. This preconditioner is also discussed in \cite{BenziGolubLiesen}.
\end{remark}



\section{Numerical experiments.}\label{sec:FE:num}
We now present some numerical experiments for our numerical space-time approaches. 
In Section~\ref{sec:num-1D} we investigate experimental orders of convergence for the conforming discretization \eqref{prob:beql_d} in one spacial dimension, in Section~\ref{ncd} we investigate experimental orders of convergence for the discretization of the saddle point problem \eqref{prob:pSaddle_disc} and test the proposed preconditioner from Section~\ref{sec:prec}.

\subsection{Conforming approximation in one spatial dimension.}
\label{sec:num-1D}
We validate the proposed conforming approach by measuring experimental convergence rates in one situation with high regularity and one situation with low regularity.
The following tests are implemented  using the Julia programming language \cite{julialang}.

We set $I=(0,1]$ and $\Omega = (0,1)$ and  consider $A=-\Delta \equiv -\partial_{xx}$. Then $V=H^1_0(\Omega)$ and $H = L^2(\Omega)$.
For the discretization we consider an equidistant grid $I_k$ with $M$ points in time and an equidistant grid $\Omega_h$ with $N$ points in space. 
Then the grid size with respect to space is $h=(N-1)^{-1}$, with respect to time is $k=(M-1)^{-1}$, and the combined grid size is $d = \sqrt{h^2+k^2}$.

In the following paragraphs we define analytical solutions to \eqref{prob:P} and calculate numerical approximations using \eqref{prob:beql_d}.
For all examples we  evaluate the errors 
\begin{align*}
\|y_d-y\|_{C^0([0,T],L^2(\Omega))}, \quad 
\|y_d-y\|_{L^2(0,T;H^1_0(\Omega))}, \quad 
\|(y_d)_t-y_t\|_{L^2(0,T;H^{-1}(\Omega))} 
\end{align*}
on a sequence of grids, where $y$ denotes the given analytical solution and $y_d$ its numerical approximation. 

Note that the last two norms together give the error contributions in the $W(0,T)$-norm and we thus expect equal convergence rates for both terms, if $y$ satisfies precisely the regularity properties of $W(0,T)$. 
We estimate rates of convergence with respect to the three grid sizes $k$, $h$, and $d$.

\paragraph{A analytical example with high regularity.}
The first example uses a smooth analytical solution and thus investigates the maximum possible convergence rates. Following \cite{FuehrerKarkulik-SpaceTimeparabolic}, we define the analytical solution
\begin{align}
   y^{(1)}(t,x) &= \sin(\pi x)\cos(\pi t)\label{eq:num:analytic-1}, 
\end{align}
so that the initial condition and the right-hand side are given by
\[
    y^{(1)}_0(x)= \sin(\pi x) \quad \text{ and } \quad
    f^{(1)}(t,x) =  \pi^2\sin(\pi x)\cos(\pi t) - \pi \sin(\pi x)\sin(\pi t).
\]
In Figure~\ref{fig:num:rates-smooth} we show the behavior of the investigated error norms for $h,k,d\to 0$, where for $d\to0$ we use equal spacing in time and space, for $h\to 0$ we fix $k = 3\cdot10^{-5}$ and for $k\to 0$ we fix $h = 4.89\cdot 10^{-4}$.
We measure 
\begin{align*}
    \|y_d-y\|_{C^0([0,T],L^2(\Omega))} \sim h^{2} + k^2, \quad
    \|y_d-y\|_{L^2(0,T;H^1_0(\Omega))} \sim h + k^2, \quad
    \|(y_d)_t-y_t\|_{L^2(0,T;H^{-1}(\Omega))} \sim h^2+k.
\end{align*}
Especially, we obtain
$\|y_d-y\|_{L^2(0,T;H^1(\Omega))} \sim  \|(y_d)_t-y_t\|_{L^2(0,T;H^{-1}(\Omega))}$ for $d =\sqrt{h^2+k^2} \to 0$, i.e.~both terms involved in the $W(0,T)$-norm asymptotically show the same convergence behavior with respect to $d\to 0$. Noting the definition of $d$, the appearing first order rate with respect to $d \to 0$ for $L^2(0,T;H^1_0(\Omega))$ and $L^2(0,T;H^{-1}(\Omega))$ is proven in \cite[Thm.~7]{FuehrerKarkulik-SpaceTimeparabolic}.

\begin{figure}
    \centering
        \begin{tikzpicture} 
    
    \pgfplotstableread[col sep = &]{imgs/ex0.equal.rates.tab}\tabSmoothEqual;

    \begin{axis}
    [
    width=0.28\textwidth,
    height=4cm,
    xlabel={$d$},
    ylabel={ },
    xmin = 0.003,
    xmax = 0.5,
    ymin = 1e-4,
    ymax = 1.0,
    ymode =log,
    xmode=log,
    legend pos = outer north east,
    ]

    
    \addplot[solid,mark=triangle*,mark size=3pt] table [x=4, y=5] from \tabSmoothEqual;
    
    \addplot[loosely dashed,mark=diamond*,mark size=3pt] table [x=4, y=9] from \tabSmoothEqual;
    
    \addplot[loosely dashed,mark=*,mark size=2pt] table [x=4, y=13] from \tabSmoothEqual;

      \addplot[loosely dashed,samples =  3, domain=0.002:1]{0.025*(x)};
      
    \addplot[solid,samples =  3, domain=0.002:1]{0.025*(x*x)};

    \end{axis}

    \end{tikzpicture}
    \begin{tikzpicture}
        \pgfplotstableread[col sep = &]{imgs/ex0.fixh.rates.tab}\tabSmoothFixh;
\begin{axis}
[
    width=0.28\textwidth,
    height=4cm,
    xlabel={$k$},
    ylabel={ },
    xmin = 0.003,
    xmax = 0.5,
    ymin = 1e-4,
    ymax = 0.1,
    ymode =log,
    xmode=log,
    legend pos = outer north east,
    ]

  
    \addplot[solid,mark=triangle*,mark size=3pt] table [x=3, y=5] from \tabSmoothFixh;
    
    \addplot[solid,mark=diamond*,mark size=3pt] table [x=3, y=9] from \tabSmoothFixh;
    
    \addplot[loosely dashed,mark=*,mark size=2pt] table [x=3, y=13] from \tabSmoothFixh;

      \addplot[loosely dashed,samples =  3, domain=0.002:1]{0.25*(x)};
      
    \addplot[solid,samples =  3, domain=0.002:1]{0.095*(x*x)};

\end{axis}
    \end{tikzpicture}
    \begin{tikzpicture}
\pgfplotstableread[col sep = &]{imgs/ex0.fixk.rates.tab}\tabSmoothFixk;
\begin{axis}
[
    width=0.28\textwidth,
    height=4cm,
    xlabel={$h$},
    ylabel={ },
    xmin = 0.003,
    xmax = 0.5,
    ymin = 1e-4,
    ymax = 1.0,
    ymode =log,
    xmode=log,
    legend pos = outer north east,
    ]

  
    \addplot[solid,forget plot] table [x=2, y=5] from \tabSmoothFixk;
    \addplot[only marks,mark=triangle*,mark size=3pt] table [x=2, y=5] from \tabSmoothFixk;
     \addlegendentry{$\eta_{C([0,T],L^2(\Omega))}$};
    
    \addplot[loosely dashed,forget plot] table [x=2, y=9] from \tabSmoothFixk;
    \addplot[only marks,mark=diamond*,mark size=3pt] table [x=2, y=9] from \tabSmoothFixk;
   \addlegendentry{$\eta_{L^2(0,T;H^1_0(\Omega))}$};
    
    \addplot[solid,forget plot] table [x=2, y=13] from \tabSmoothFixk;
    \addplot[only marks,mark=*,mark size=2pt] table [x=2, y=13] from \tabSmoothFixk;
 \addlegendentry{$\eta_{L^2(0,T;H^{-1}(\Omega))}$};

      \addplot[loosely dashed,samples =  3, domain=0.002:1]{0.185*(x)};
         \addlegendentry{Order 1}
         
    \addplot[solid,samples =  3, domain=0.002:1]{0.025*(x*x)};
  \addlegendentry{Order 2}
\end{axis}
    \end{tikzpicture}

    \caption{Error plots for the smooth example \eqref{eq:num:analytic-1}.
    Here $\eta_{C([0,T],L^2(\Omega))} = \|y-y_d\|_{C([0,T],L^2(\Omega))}$,
    $\eta_{L^2(0,T;H^1_0(\Omega))} = \|y-y_d\|_{L^2(0,T;H^1_0(\Omega))}$, and
    $\eta_{L^2(0,T;H^{-1}(\Omega))} = \|y_t-(y_d)_t\|_{L^2(0,T;H^{-1}(\Omega))}$. 
    For the left plot we use $h=k$, for the middle plot we fix $h=4.89\cdot10^{-4}$ and for the right plot we fix  $k = 3\cdot10^{-5}$. We observe 
    $ \|y-y_d\|_{C([0,T],L^2(\Omega))} \sim h^2 + k^2$,
    $\|y-y_d\|_{L^2(0,T;H^1_0(\Omega))} \sim h+k^2$, and 
    $\|y_t-(y_d)_t\|_{L^2(0,T;H^{-1}(\Omega))} \sim h^2 + k$.
    }
    \label{fig:num:rates-smooth}
\end{figure}


\paragraph{An analytical example with low regularity.}

The second example has an analytical solution with low regularity. Here we use $f \in L^\infty(0,T;\mathcal  (C(\overline \Omega))^\star)$.
Since $\Omega \subset \RR$ we have $H^1_0(\Omega) \hookrightarrow C(\overline \Omega)$ and thus $(C(\overline \Omega))^\star \hookrightarrow H^{-1}(\Omega)$.
We set 
\begin{align}
   y^{(2)}(t,x) &= (0.5-|x-0.5|)(|t-0.5|+0.5)\label{eq:num:analytic-2},
\end{align}
which is the analytical solution of \eqref{eq:parab} related to the data
\begin{align*}
    y^{(2)}_0(x) &= (0.5-|x-0.5|), 
    & f^{(2)}(t,x)   &= (0.5-|x-0.5|)\sign(t-0.5) + 2 \delta_{0.5}(x)(|t-0.5|+0.5),
\end{align*}
where $\delta_{0.5}(x)$ denotes the Dirac point measure concentrated at $x=0.5$. To prevent from superconvergence, we use meshes, that do not contain the point $x=0.5$ as discretization point.

In Figure~\ref{fig:num:rates-rough} we show the behavior of the investigated error norms for $h,k,d\to 0$, where for $d \to 0$ we consider equal spacing for time and space, for $k\to 0$ we fix $h=4.89\cdot 10^{-4}$ and for $h \to 0$ we fix $k=7.63\cdot 10^{-6}$.
We obtain 
\begin{align*}
   &\|y_d-y\|_{L^2(0,T;H^1_0(\Omega))} \sim h^{1/2}+k^p, \text{ where }p\geq 1/2 ,  \quad \|(y_d)_t-y_t\|_{L^2(0,T;H^{-1}(\Omega))} \sim h^2 + k^{1/2},\\
     &\|y_d-y\|_{C^0([0,T],L^2(\Omega))} \sim h^{3/2} + k.
\end{align*}
Especially, we again obtain
$\|y_d-y\|_{L^2(0,T;H^1_0(\Omega))} \sim  \|(y_d)_t-y_t\|_{L^2(0,T;H^{-1}(\Omega))}$ for $d \to 0$.

\begin{figure}
    \centering
        \begin{tikzpicture} 
    
    \pgfplotstableread[col sep = &]{imgs/ex1.equal.rates.tab}\tabRoughEqual;

    \begin{axis}
    [
    width=0.28\textwidth,
    height=4cm,
    xlabel={$d$},
    ylabel={ },
    xmin = 0.003,
    xmax = 0.5,
    ymin = 1e-4,
    ymax = 1.0,
    ymode =log,
    xmode=log,
    legend pos = outer north east,
    ]

    
    \addplot[loosely dashed,mark=triangle*,mark size=3pt] table [x=4, y=5] from \tabRoughEqual;
    
    \addplot[dotted,mark=diamond*,mark size=3pt] table [x=4, y=9] from \tabRoughEqual;
    
    \addplot[dotted,mark=*,mark size=2pt] table [x=4, y=13] from \tabRoughEqual;

      \addplot[loosely dashed,samples =  3, domain=0.002:1]{0.025*(x)};
      
    \addplot[dotted,samples =  3, domain=0.002:1]{0.095*(x^(0.5))};
    

    \end{axis}

    \end{tikzpicture}
    \begin{tikzpicture}
        \pgfplotstableread[col sep = &]{imgs/ex1.fixh.rates.tab}\tabRoughFixh;
\begin{axis}
[
    width=0.28\textwidth,
    height=4cm,
    xlabel={$k$},
    ylabel={ },
    xmin = 0.003,
    xmax = 0.5,
    ymin = 3e-4,
    ymax = 1e-1,
    ymode =log,
    xmode=log,
    legend pos = outer north east,
    ]

  
    \addplot[loosely dashed,mark=triangle*,mark size=3pt] table [x=3, y=5] from \tabRoughFixh;
    
    \addplot[loosely dashed,mark=diamond*,mark size=3pt] table [x=3, y=9] from \tabRoughFixh;
    
    \addplot[dotted,mark=*,mark size=2pt] table [x=3, y=13] from \tabRoughFixh;

      \addplot[loosely dashed,samples =  3, domain=0.002:1]{0.30*(x)};
      
    \addplot[dotted,samples =  3, domain=0.002:1]{0.01*(x^(0.5))};
    

\end{axis}
    \end{tikzpicture}
    \begin{tikzpicture}
\pgfplotstableread[col sep = &]{imgs/ex1.fixk.rates.tab}\tabRoughFixk;
\begin{axis}
[
    width=0.28\textwidth,
    height=4cm,
    xlabel={$h$},
    ylabel={ },
    xmin = 0.008,
    xmax = 0.5,
    ymin = 1e-5,
    ymax = 1.0,
    ymode =log,
    xmode=log,
    legend pos = outer north east,
    ]

  
    \addplot[dashed,forget plot] table [x=2, y=5] from \tabRoughFixk;
    \addplot[only marks,mark=triangle*,mark size=3pt] table [x=2, y=5] from \tabRoughFixk;
  \addlegendentry{$\eta_{C([0,T],L^2(\Omega))}$};
    
    \addplot[dotted,forget plot] table [x=2, y=9] from \tabRoughFixk;
    \addplot[only marks,mark=diamond*,mark size=3pt] table [x=2, y=9] from \tabRoughFixk;
    \addlegendentry{$\eta_{L^2(0,T;H^1_0(\Omega))}$};
    
    \addplot[solid,forget plot] table [x=2, y=13] from \tabRoughFixk;
    \addplot[only marks,mark=*,mark size=2pt] table [x=2, y=13] from \tabRoughFixk;
     \addlegendentry{$\eta_{L^2(0,T;H^{-1}(\Omega))}$};

    \addplot[dotted,samples =  3, domain=0.001:1]{0.2*(x^(0.5))};
          \addlegendentry{Order 0.5}
    
    \addplot[loosely dashed,samples=3,domain=0.001:1]{10^(-6)*x};
    \addlegendentry{Order 1}
    
   \addplot[ dashed,samples=3,domain=0.001:1]{0.115*x^(1.5)};
    \addlegendentry{Order 1.5}
   
        \addplot[solid,samples =  3, domain=0.001:1]{0.005*(x*x)};
        \addlegendentry{Order 2}

\end{axis}
    \end{tikzpicture}
    \caption{Error plots for the non-smooth example \eqref{eq:num:analytic-2}.
    Here $\eta_{C([0,T],L^2(\Omega))} = \|y-y_d\|_{C([0,T],L^2(\Omega))}$,
    $\eta_{L^2(0,T;H^1_0(\Omega))} = \|y-y_d\|_{L^2(0,T;H^1_0(\Omega))}$, and
    $\eta_{L^2(0,T;H^{-1}(\Omega))} = \|y_t-(y_d)_t\|_{L^2(0,T;H^{-1}(\Omega))}$. 
    For the left plot we use $h=k$, for the middle plot we fix  $h=4.89\cdot10^{-4}$ and for the right plot we fix $k = 7.63\cdot10^{-6}$. We observe 
    $ \|y-y_d\|_{C([0,T],L^2(\Omega))} \sim h^{3/2} + k$,
    $\|y-y_d\|_{L^2(0,T;H^1_0(\Omega))} \sim h^{1/2}+k^p$, where $p\geq 1/2$ and 
    $\|y_t-(y_d)_t\|_{L^2(0,T;H^{-1}(\Omega))} \sim h^2 + k^{1/2}$. We are not able to fix a smaller $h$ to evaluate the error with respect to $k$ for $\|y-y_d\|_{L^2(0,T;H^1_0(\Omega))}$ due to technical limitations, but the minimum rate $1/2$ can be seen from the results for $d\to 0$.}
    \label{fig:num:rates-rough}
\end{figure}

\FloatBarrier
\subsection{Nonconforming discretization in two spatial dimensions.}
\label{ncd}
The implementation in two spatial dimensions is done in \texttt{Python}.
We use the finite element package \texttt{dolfinx} \cite{BarattaEtal2023,ScroggsEtal2022,BasixJoss} together with \texttt{UFL} \cite{AlnaesEtal2014}. The linear algebra backend is provided by \texttt{PETSc} \cite{petsc1,petsc2}.

\begin{figure}
    \centering
    \includegraphics[width=0.16\textwidth, trim={27pt 200pt 27pt 200pt},clip]{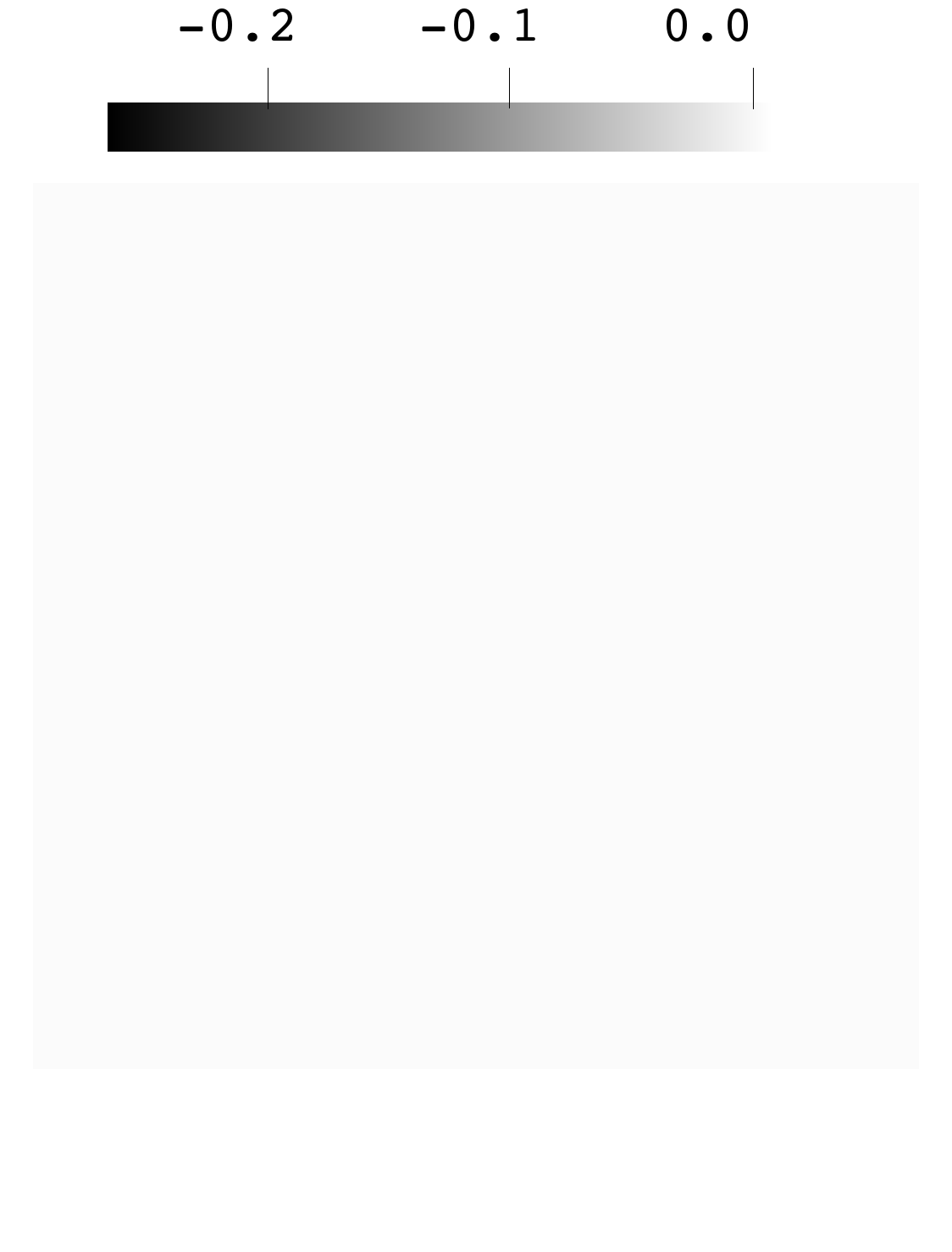}
    \includegraphics[width=0.16\textwidth, trim={27pt 200pt 27pt 200pt},clip]{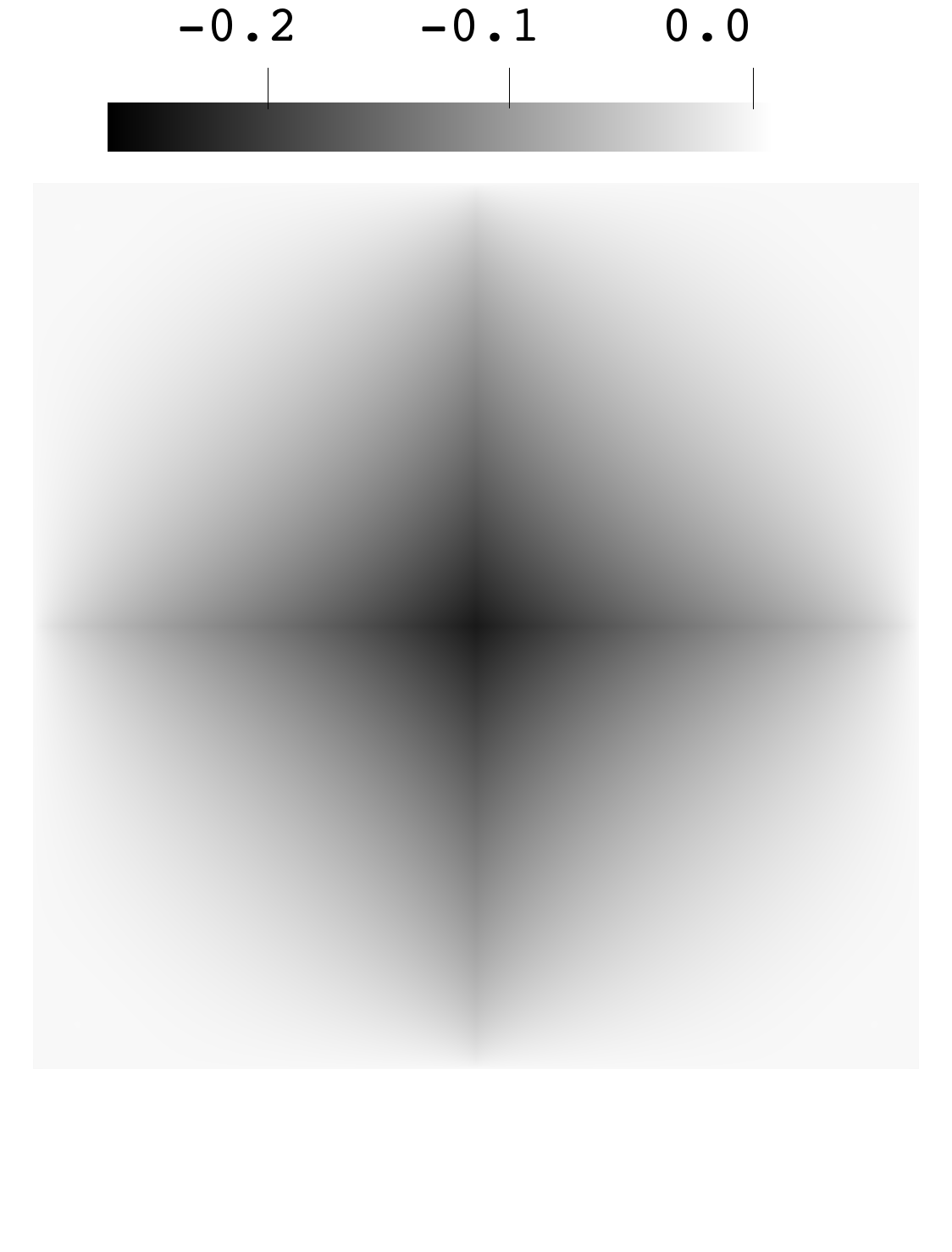}
    \includegraphics[width=0.16\textwidth, trim={27pt 200pt 27pt 200pt},clip]{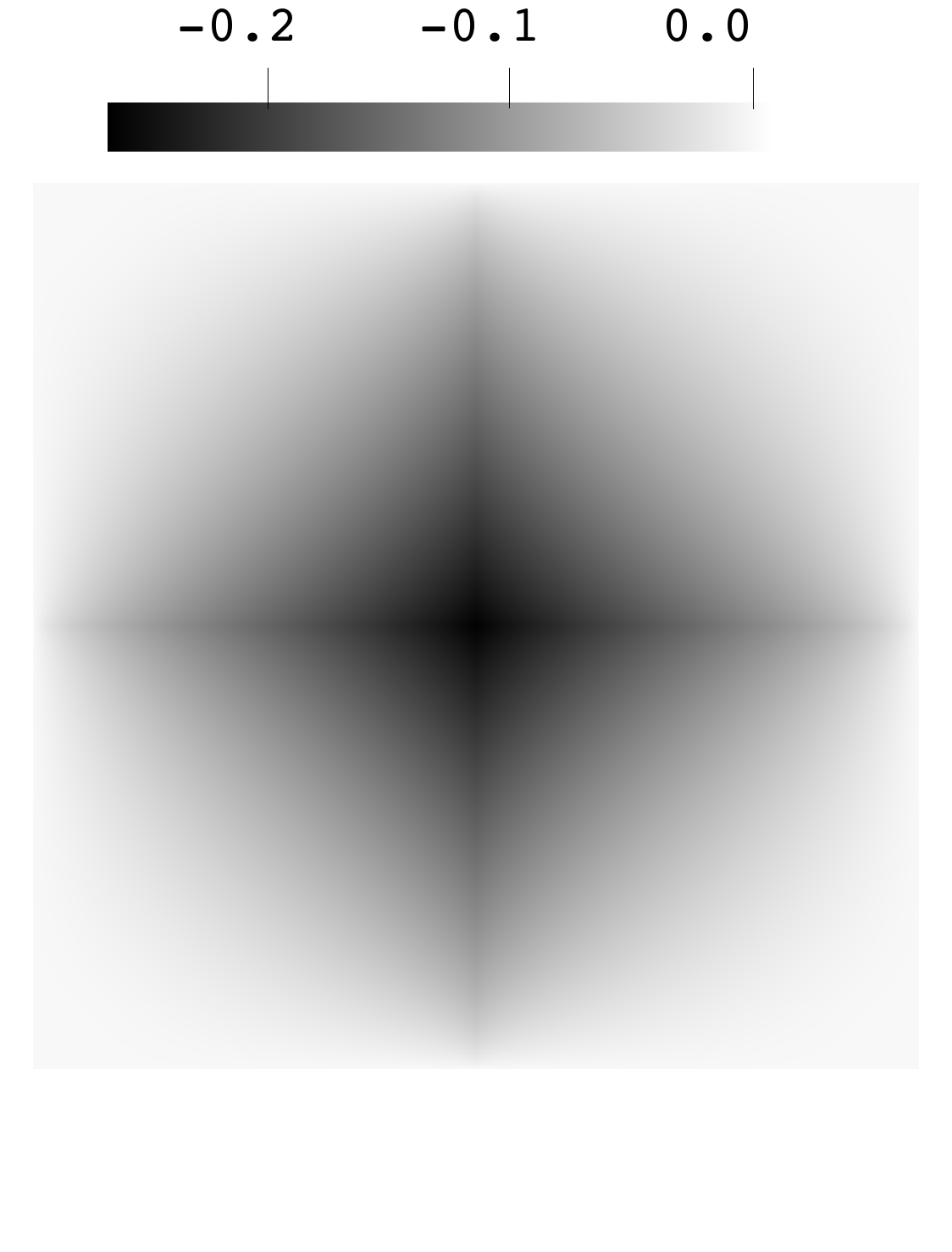}
    \includegraphics[width=0.16\textwidth, trim={27pt 200pt 27pt 200pt},clip]{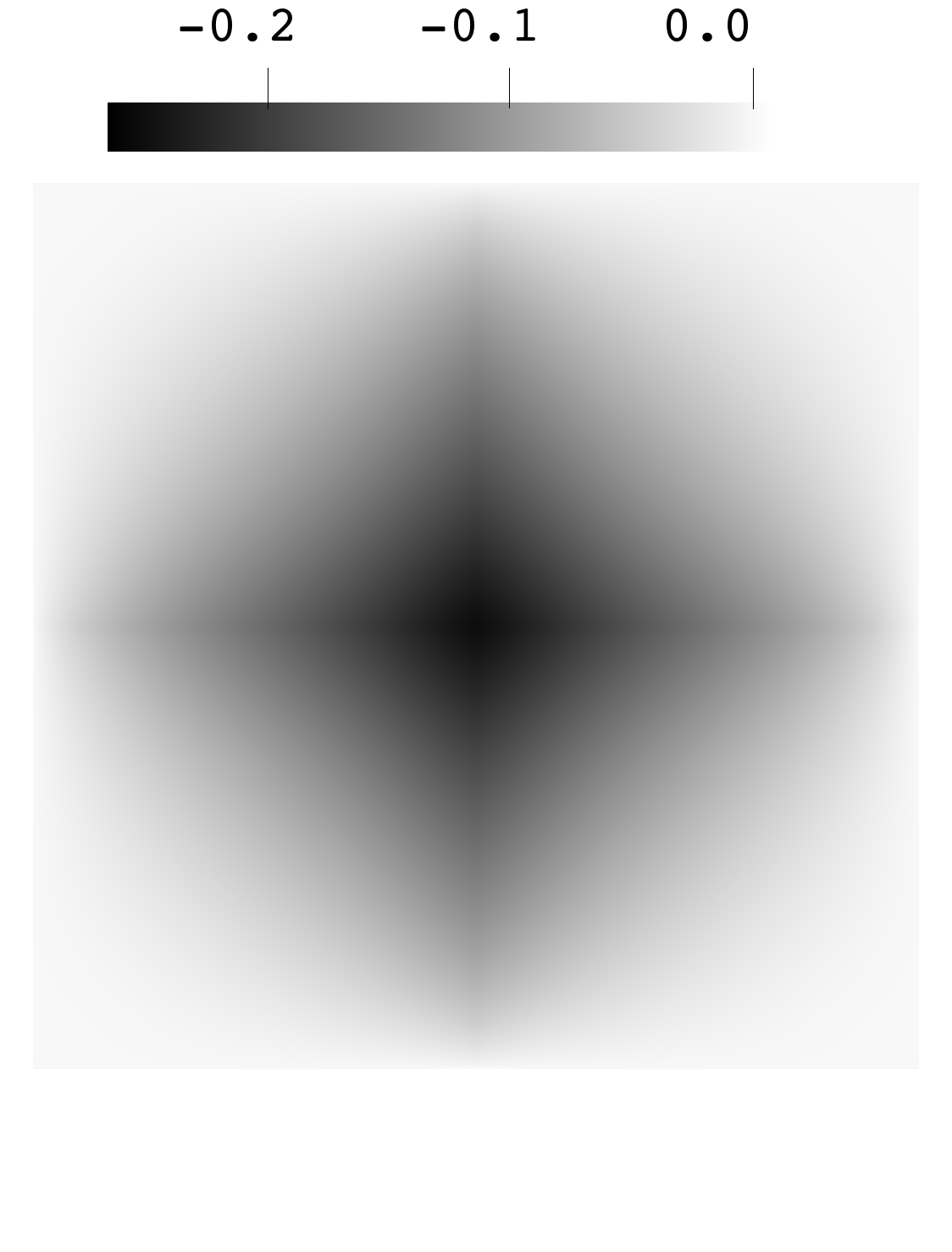}
    \includegraphics[width=0.16\textwidth, trim={27pt 200pt 27pt 200pt},clip]{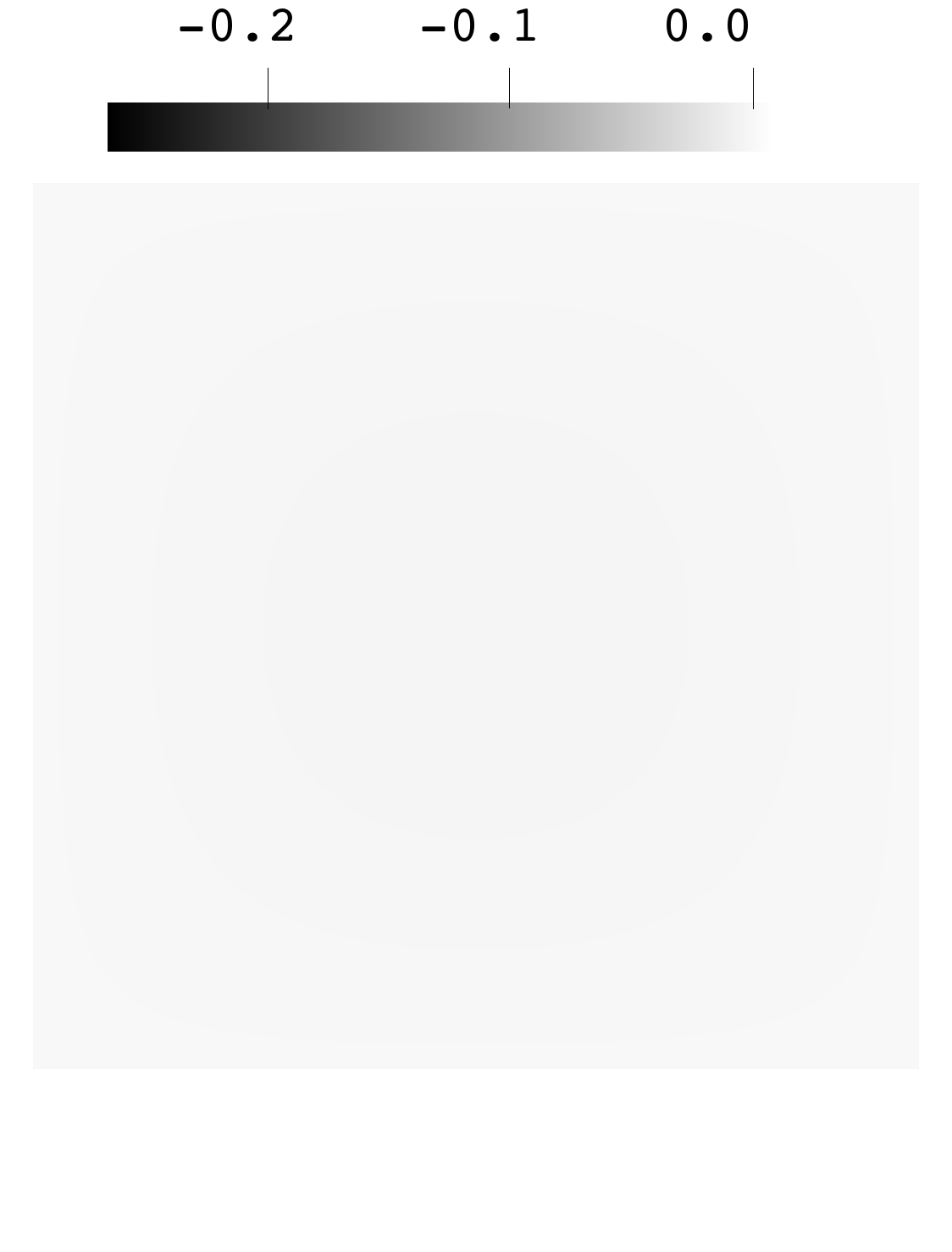}
    \includegraphics[width=0.16\textwidth, trim={27pt 200pt 27pt 200pt},clip]{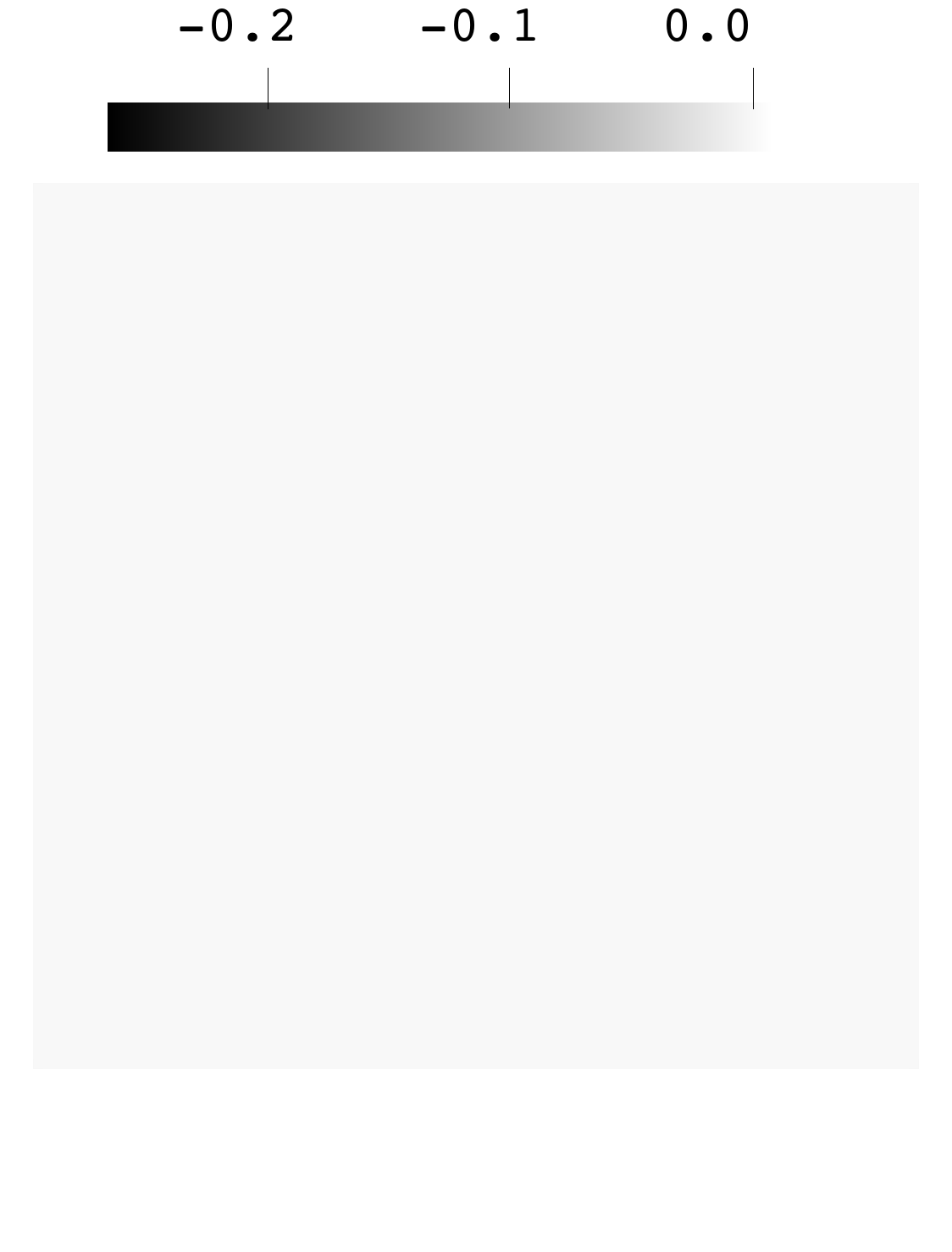}
    \caption{Plot of the numerical approximation of the  solution of example \eqref{eq:bsprechteseite0}--\eqref{eq:bsprechteseite2} at times $t = 0.0$, $0.1$, $0.3$, $0.5$, $0.7$ and $0.9$ (left to right). Black indicates negative values. The solution shows kinks at $(x_1,0.5)$ and $(0.5,x_2)$ and rapidly decays to zero after $t=0.5$ with a kink at $t=0.5$. This can be seen in more detail in Figure \ref{fig:timeevoly}.}
    \label{fig:plot_true_sol}
\end{figure}

We investigate the proposed preconditioner and the convergence rate at hand with some non-smooth example in two spatial dimensions. We set $\Omega = (0,1)^2$, $I = (0, T] = (0,1]$ and $A = - \Delta$. 
We introduce an equidistant discretization of $I$ with $M$ vertices as well as an equidistant regular triangulation of $\Omega$ with $M$ vertices in each spatial direction and $N=M^2$ total vertices.
We choose $y_0 = 0$ as well as a separable source term $f(t,x) = f^t(t) f^x(x) \in L^\infty(0,T;(W^{1,1}(\Omega))^\star)$,
with 
\begin{align}\label{eq:bsprechteseite0}
    f^t(t) := \begin{cases}1, & \text{if } t \leq 0.5 \\ 0, & \text{else}\end{cases}
\end{align}
and for arbitrary $v \in V$ we choose
\begin{align}\label{eq:bsprechteseite1}
    \langle f^x, v \rangle_{V^\star, V} = \int_\Omega f_0 v + (f_1, f_2)^T \nabla v \dx,
\end{align}
where $f_0 \equiv 1$ and for $i=1,2$
\begin{align}\label{eq:bsprechteseite2}
    f_i(x) := 
    \begin{cases}
    1, & \text{if } x_i \geq 0.5, \\ 
    0, & \text{else}.
    \end{cases} 
\end{align}
In Figure~\ref{fig:plot_true_sol} we show a numerical approximation to $y$ using this data $(y_0,f)$.

\begin{figure}
    \centering
    \begin{tikzpicture} 
\pgfplotstableread[col sep = comma]{sim/time_evol.csv}\tabTimeEvol;
\begin{axis}
[
width=0.49\textwidth,
height=4cm,
xlabel={$t$},
legend pos = outer north east,
]

\addplot[thick] table [x index = 0, y index = 1] from \tabTimeEvol;

\addlegendentry{$y(t,x)|_{x=(0.5,0.5)}$}

\end{axis}
\end{tikzpicture}
    \caption{Evolution over time of $y(t,x)$ at $x = (0.5,0.5)$. We observe a kink at $t = 0.5$, which is expected by \eqref{eq:bsprechteseite0}.}
    \label{fig:timeevoly}
\end{figure}

\paragraph{Investigation of the experimental order of convergence.}
To calculate experimental orders of convergence we calculate $y_d$ on a sequence of grids with $M$ vertices in time and $N=M^2$ vertices with respect to space. To prevent from super-convergence  we choose $M$ as even number such that the discontinuities in $f$ are not resolved. 

An analytical solution $y$ is not available. We therefore approximate $y$ by a numerical solution with high resolution.
For this we additionally solve the problem for an odd value of $M$, chosen as twice the finest level minus one, such that the discontinuities are resolved, and consider this as approximation of the analytical solution $y$. We plot $y$ in Figure \ref{fig:plot_true_sol} and $y$ evaluated at $x_1 = x_2 = 0.5$ over time in Figure \ref{fig:timeevoly}.

\begin{figure}
    \centering
    \begin{tikzpicture} 
\pgfplotstableread[col sep = comma]{sim/results_conv.csv}\tabConv;
\begin{axis}
[
width=0.49\textwidth,
height=4cm,
xlabel={$d$},
ylabel={},
xmin=1e-2,
xmax=1.5e-1,
ymode = log,
xmode = log,
legend pos = outer north east,
]
\addplot[loosely dashed,samples =  3, domain=0.01:0.05]{10^(-0.8432+0.5)*x};
\addlegendentry{$\mathcal{O}(d)$}

\addplot[loosely dashed,mark=diamond*,mark size=3pt,mark repeat=6,mark phase=5] table [x index = 2, y index = 4] from \tabConv;
\addlegendentry{$\vert \vert y_d - y \vert \vert_{L^2(0,T; H^1_0(\Omega))}$}

\addplot[loosely dashed, mark=*,mark size=2pt, mark phase=5, mark repeat =6] table [x index = 2, y index = 5] from \tabConv;
\addlegendentry{$\vert \vert (y_d)_t - y_t \vert \vert_{L^2(0,T; H^{-1}(\Omega))}$}


\addplot[dashed, mark=triangle*, mark size=3pt, mark phase=5, mark repeat =6] table [x index = 2, y index = 3] from \tabConv;
\addlegendentry{$\vert \vert y_d - y \vert \vert_{C([0,T], L^2(\Omega))}$}


\addplot[dashed,samples =  3, domain=0.01:0.05]{10^(-1.9709-0.5)*x^1.5};
\addlegendentry{$\mathcal{O}(d^{1.5})$}

\end{axis}
\end{tikzpicture}
    \caption{Finite element errors for $d \rightarrow 0$ in the given example for two spatial dimensions. For fine enough $d$ we obtain a convergence of order of 1.5 with respect to $d$ in the  $C([0,T],L^2(\Omega))$ norm and a convergence order of 1.0 with respect to $d$ in the $L^2(0,T;H^1_0(\Omega))$ and $L^2(0,T;H^{-1}(\Omega))$ norm.}
    \label{fig:num2D:EOC}
\end{figure}

In Figure \ref{fig:num2D:EOC} we shown the behavior of the investigated error norm for $d \rightarrow 0$. We obtain
\begin{align}
    \vert \vert y_d - y \vert \vert_{L^2(0,T;H^1_0(\Omega))} \sim d, \quad \vert \vert (y_d)_t - y_t \vert \vert_{L^2(0,T;H^{-1}(\Omega))} \sim d, \quad \vert \vert y_d - y \vert \vert_{C([0,T],L^2(\Omega))} \sim d^{1.5}.
\end{align}
Especially we observe $\vert \vert y_d - y \vert \vert_{L^2(0,T;H^1_0(\Omega))} \sim \vert \vert (y_d)_t - y_t \vert \vert_{L^2(0,T;H^{-1}(\Omega)} \sim d$, which can be expected, since these two norms reproduce the error in the $W(0,T)$ norm and our example satisfies precisely the regularity properties of $W(0,T)$.

\paragraph{Investigation of the proposed preconditioner.}
We investigate the proposed preconditioners in combination with GMRES \cite{gmres}, LGMRES \cite{lgmres} and MINRES \cite{minres}, i.e.~we use the triagonal preconditioner \eqref{eq:P_left} for GMRES and LGMRES and the symmetric and positive definite preconditioner \eqref{eq:P_sym} for MINRES. 
Here, we use code fragments from the corresponding \texttt{SciPy} implementations \cite{2020SciPy-NMeth} and \cite{Golub}. 
The preconditioner is implemented using the \texttt{LinearOperator} interface from \texttt{SciPy} as well as \texttt{PETSc} \cites{petsc1}{petsc2}.
We use the algebraic multigrid method \cite{STUBEN1983419} provided by \texttt{PETSc} with ten Chebyshev smoothing steps \cite{houwen}.

For all methods we use a relative tolerance of $10^{-5}$  
and allow at most $10^6$ iterations. The convergence is checked with the unpreconditioned residual to allow a meaningful comparison of the methods with and without preconditioner. 
For GMRES we consider a restart every 30 iterations.
For LGMRES we consider the restart every 30 iterations carrying 3 augmentation vectors.

\begin{figure}
    \centering
    \begin{tikzpicture} 
\pgfplotstableread[col sep = comma]{sim/results_prec.csv}\tabPrec;
\begin{groupplot}[
    group style={
        group size=3 by 1,
        horizontal sep=8pt,
        ylabels at=edge left,
    },
    footnotesize,
    width=0.35\textwidth,
    height=5cm,
    xlabel=\# DOFs,
    max space between ticks=20pt,
    ymin=0.2,
    ymax=10^6,
    ymode=log,
    xmode=log,
]

\nextgroupplot[title=GMRES, ylabel=\# Iterations]
\addplot[loosely dashed, mark=*,mark size=2pt, mark phase=92, mark repeat =999] table [x index = 0, y index = 3] from \tabPrec;
\addplot[loosely dashed,mark=diamond*,mark size=3pt,mark phase=92, mark repeat =999] table [x index = 0, y index = 4] from \tabPrec;
\addplot[only marks, mark=*,mark size=2pt] coordinates {(501.54243680330615,22)(11962.586877716014,7646)};
\addplot[only marks, mark=diamond*,mark size=3pt] coordinates {(501.54243680330615,6)(11962.586877716014,8)(343853.62352610176,8)};

\nextgroupplot[title=LGMRES, yticklabels={}, legend style={at={($(0,0)+(0.12\textwidth,3cm)$)},legend columns=4,fill=none,draw=black,anchor=center,align=center, yshift=1.5cm}]
\addplot[loosely dashed, mark=*,mark size=2pt, mark phase=999, mark repeat =999] table [x index = 0, y index = 5] from \tabPrec;
\addlegendentry{Without preconditioner \,}
\addplot[loosely dashed,mark=diamond*,mark size=3pt,mark phase=92, mark repeat =999] table [x index = 0, y index = 6] from \tabPrec;
\addlegendentry{With preconditioner \,}
\addplot[only marks, mark=*,mark size=2pt] coordinates {(501.54243680330615,4.878653791974776)(11962.586877716014,36.08278079727419)(343853.62352610176,266.8701501643408)};
\addplot[only marks, mark=diamond*,mark size=3pt] coordinates {(501.54243680330615,2)(11962.586877716014,2)(343853.62352610176,2)};

\nextgroupplot[title=MINRES, yticklabels={}]
\addplot[loosely dashed, mark=*,mark size=2pt, mark phase=999, mark repeat =999] table [x index = 0, y index = 7] from \tabPrec;
\addplot[loosely dashed,mark=diamond*,mark size=3pt,mark phase=92, mark repeat =999] table [x index = 0, y index = 8] from \tabPrec;
\addplot[only marks, mark=*,mark size=2pt] coordinates {(501.54243680330615,55)(11962.586877716014,420)(343853.62352610176,1666)};
\addplot[only marks, mark=diamond*,mark size=3pt] coordinates {(501.54243680330615,9)(11962.586877716014,13)(343853.62352610176,11)};

\end{groupplot}
\end{tikzpicture}
    \caption{Comparison of the number of iterations of the given iterative methods with and without preconditioner. If the preconditioner is applied, the number of iterations in all three compared methods remains approximately constant as desired. Without preconditioner, the LGMRES and the MINRES method behave rather similar, whereas the standard GMRES fails to converge within $10^6$ maximal iterations for large $M$.}\label{fig:preconditioner_plot}
\end{figure}

As before, we choose a sequence of even numbers for $M$ and set $N=M^2$ leading to $\mathcal{O}(M^3)$ unknowns of the linear systems. Due to the equal spacing in time and each spatial dimension we obtain $d=\sqrt{3}(M-1)^{-1}$.
The number of iterations in dependence on the degrees of freedom are shown in Figure \ref{fig:preconditioner_plot}. 
We observe that the preconditioner in all three cases strongly reduces the needed iterations. The number of iterations stays approximately constant, whereas the methods without preconditioner require significantly more effort. Here, LGMRES and MINRES behave rather similar and the restarted GMRES fails to converge within the given relative tolerance and $10^6$ maximal iterations for large systems.


\section{Conclusion.}
In this work we introduce a least squares formulation for the solution of parabolic problems in the $L^2(0,T;V^\star) \times H$ norm requiring only minimal regularity of the data. For the numerical approximation we propose a conforming and a nonconforming Galerkin discretization and prove convergence for both approaches. In the conforming approach we rely on exact evaluations of the Riesz representers of finite element basis functions, which practically is implementable for spatially one dimensional problems. For higher space dimensions we propose a nonconforming discretization which is treated with a saddle point formulation. For the resulting discrete space-time systems we propose quasi-optimal preconditioners. Numerical experiments are presented which confirm our analytical findings.

\paragraph{Author contributions (CRediT taxonomy)}
\begin{itemize}
    \item M.H.: Idea, Project administration, Conceptualization,  Methodology, Formal analysis, Writing – original draft, Writing – review \& editing, Funding acquisition;
    \item C.K.: Project administration, Conceptualization, Methodology, Formal analysis, Software, Investigation, Writing – original draft, Writing – review \& editing, Funding acquisition;
    \item M.S.: Methodology, Formal analysis, Software, Investigation, Data curation, Visualization, Writing – review \& editing.
\end{itemize}

\paragraph{Funding.}

The first author acknowledges funding of the project 
\textit{Ein nichtglatter Phasenfeld Zugang für Formoptimierung mit instationären Fluiden} by the German Research foundation  within the Priority Programme 1962 under project number 423457678.
The first and second author also acknowledge funding of the project 
\textit{Fluiddynamische Formoptimierung mit Phasenfeldern und Lipschitz-Methoden} by the German Research foundation under project number 543959359.


\printbibliography

\end{document}